\documentclass[11pt]{article}
\usepackage{graphicx}
\usepackage{amsmath,amssymb,amsthm,amsfonts,mathrsfs}
\usepackage{amssymb}
\usepackage{hyperref}
\usepackage[a4paper]{geometry}
\newtheorem{thm}{Theorem}[section]

\newtheorem{lem}[thm]{Lemma}

\theoremstyle{definition}
\newtheorem{defn}{Definition}[section]
\theoremstyle{remark}
\newtheorem{rem}{Remark}[section]

\numberwithin{equation}{section}

\DeclareMathSymbol{\C}{\mathalpha}{AMSb}{"43}

\newcommand{\eps}{\varepsilon}

\newcommand{\alp}{\alpha}

\newcommand{\dx}{\,\mathrm{d}x}

\newcommand{\R}{{\mathbb{R}}}

\newcommand{\inte}{\int_{\mathbb{R}^N}}
\newcommand{\intB}{\int _{B_\delta (x_{\eps})}}
\newcommand{\intPB}{\int _{\partial B_\delta (x_{\eps})}}

\newcommand{\bsub}{\begin{subequations}}
\newcommand{\esub}{\end{subequations}$\!$}

\begin{document}
\title{Uniqueness of Single Peak Solutions for  Coupled Nonlinear Gross-Pitaevskii Equations with  Potentials }
\author{ Xiaoyu Zeng$$\ \ and Huan-Song Zhou$$\\
   \small \it Department of Mathematics, Wuhan University of Technology, Wuhan 430070, P.R. China \\
  }
\date{\today}

\smallbreak\maketitle
\begin {abstract}
For a couple of singularly perturbed Gross-Pitaevskii equations, we first prove that the single peak solutions, if they  concentrate on the same point, are unique provided that the Taylor's expansion of potentials around the concentration point is in the same order along all directions. Among other assumptions, our results  indicate that the peak solutions obtained in \cite{IT,MPS,WS}  are unique. Moreover, for the radially symmetric ring-shaped potential, which attains its minimum at the spheres $\Gamma_i:=\{x\in\R^N:|x|=A_i>0\},i=1,2,\cdots,l,$  and is totally degenerate  in the tangential space of  $\Gamma_i$, we prove that the positive ground state is cylindrically symmetric and is unique up to rotations around the origin. Aa far as we know, this is the first uniqueness result for ground states under radially symmetric but non-monotonic potentials.
\end {abstract}

\vskip 0.2truein
\noindent {\it Keywords:} Gross-Pitaevskii equations; single peak solutions;  ring-shaped potentials;  uniqueness.
\vskip 0.2truein
\noindent  {\em MSC(2010): 35J47, 35J50, 46N50}
\vskip 0.2truein


\section{Introduction}
In this paper, we consider the following time-independent coupled nonlinear Gross-Pitaevskii equations
\begin{equation}\label{eq1.1}
\begin{cases}
-\eps^2\Delta u_{1} +V_1(x)u_{1} =a_1u_{1}^3 +\beta u_{2}^2 u_{1}   \,\ \mbox{in}\,\  \R^N,\ 1\leq N\leq 3,\\
-\eps^2\Delta u_{2} +V_2(x)u_{2} =a_2u_{2}^3 +\beta u_{1}^2 u_{2}   \,\ \mbox{in}\,\  \R^N.\,\
\end{cases}
\end{equation}
System \eqref{eq1.1} mainly arises in the theoretical study of two-component Bose-Einstein condensation(BEC), which has been extensively studied by mathematicians and physicists, in view of the experimental realization of two component BEC  for trapped alkali atomic gases \cite{HMEWC}. Here $(u_1(x), u_2(x))$ is related to the macroscopic wave function of the ultracold atom system, $V_i(x)\in C(\R^N; \R^+)$($i=1,2,$) denotes the external trapping potential, the absolute value of the parameter $|a_i|(i=1,2,)$ describes the interaction strength among  the atoms in  each component, while $|\beta|$ denotes the interaction strength  for the atoms between different  components. Meanwhile, the signs of $a_1,a_2$ and $\beta$ denote the interactions are attractive or repulsive, respectively.

The classification on the existence of solutions of equation \eqref{eq1.1} is  an interesting but quit difficult topic, for which depends sensitively on the parameters $a_1,a_2, \beta$, as well as the potentials $V_i(x) (i=1,2,)$. Moreover, comparing with the case of single equation, another difficulty of studying \eqref{eq1.1} lies in how to distinguish the non-trivial solutions ( $u_1\not=0$ and $u_2\not=0$) from the  semi-trivial ones (i.e., one of $u_1$ and $u_2$ equals to zero). As a consequence, problem \eqref{eq1.1} has received munch attention  in the last decade. In particular, after the pioneering work of Lin and Wei \cite{LWCMP} which considered the existence of nontrivial solutions for  the following  autonomous system corresponding to \eqref{eq1.1} with $V_i(x)=\mu_i>0$ and $\eps=1$,
\begin{equation}\label{eq1.01}
\begin{cases}
-\Delta u_{1} +\mu_1u_{1} =a_1u_{1}^3 +\beta u_{2}^2 u_{1}   \,\ \mbox{in}\,\  \R^N,\ 1\leq N\leq 3,\\
-\Delta u_{2} +\mu_2u_{2} =a_2u_{2}^3 +\beta u_{1}^2 u_{2}   \,\ \mbox{in}\,\  \R^N.
\end{cases}
\end{equation}
For instance, Sirakov \cite{Si} obtained the ground states of \eqref{eq1.01} by searching for minimizers of the energy functional corresponding to \eqref{eq1.01} over the Nehari manifold. In particular, he proved that, there exist $\beta_1=\beta_1(a_1,a_2,\mu_1,\mu_2)\in(0,\min\{a_1,a_2\}]$  and $\beta_2=\beta_2(a_1,a_2,\mu_1,\mu_2)\in[\max\{a_1,a_2\},+\infty)$ (see \cite[Theorem 1.2]{Si} for the explicit form of $\beta_i$ ) such that, \eqref{eq1.01} has a positive ground state provided that $\beta\in (0,\beta_1)\cup(\beta_2,+\infty)$, while \eqref{eq1.01} has no positive solution provided that $\beta\in [\min\{a_1,a_2\}, \max\{a_1,a_2\}]$. Moreover, for the  special case of $\mu_1=\mu_2=\mu>0$, i.e.,  for the following system
\begin{equation}\label{sys1}
\arraycolsep=1.5pt
\left\{\begin{array}{lll}
&-\Delta u +\mu u=a_1 u ^3 + \beta  v ^2 u   \,\ \mbox{in}\,\  \R^N,\\ [2.0mm]
& -\Delta v +\mu v =a_2 v ^3 +  \beta u ^2 v   \,\ \mbox{in}\,\  \R^N,\,\\
&u,v\in H^1(\R^N),
\end{array}\right.
\end{equation}
 he showed  that $\beta_1=\min\{a_1,a_2\}$ and $\beta_2=\max\{a_1,a_2\}$, and  $(w_1(x),w_2(x))$ with
\begin{equation}\label{sys2}
w_i(x)=\sqrt{\gamma_i\mu}w(\sqrt{\mu}x)\ \text{ where }\gamma_i=\frac{\beta-a_j}{\beta^2-a_1a_2}\text{ and }(i,j)=(1,2) \text{ or }(2,1),
\end{equation}
is the positive ground state of \eqref{sys1}. Here  $w(|x|)$ is the unique radially symmetric positive solution of the field equation \cite{Kwong}
\begin{equation}\label{equ:w}
\Delta w-w+w^3=0,\,\,\  w \in H^1(\R^N).
\end{equation}
However, as pointed out in \cite{Si}, when $\mu_1\not=\mu_2$,  $\beta_1$ and $\beta_2$ are not the optimal endpoints which determine the existence and non-existence of nonnegative nontrivial solutions of \eqref{eq1.01}. To find out
 the optimal ranges for the existence of nonnegative nontrivial solutions is a quite interesting problem. Recently, the authors of \cite{BZZ,WZZ}  make some progress on this aspect, and obtained some new intervals  about $\beta$ for the existence and nonexistence of nonnegative solutions.

 Sirakov conjectured in \cite{Si} that, up to translations, $(w_1(x),w_2(x))$  with the form of \eqref{sys2}  is the unique positive solution of \eqref{sys1} when \begin{equation}\label{eq1.06}\beta\in(0,\min\{a_1,a_2\})\cup(\max\{a_1,a_2\},+\infty).\end{equation}
Wei and Yao \cite{WY} affirmed  this conjecture for  $
\beta\in(\max\{a_1,a_2\},+\infty)$, or $\beta>0$ and small enough.  Chen  and Zou \cite{CZ2} proved this conjecture for the case   that  $\beta>0$ is less  but close to $\min\{a_1,a_2\}$.  Afterwards, they considered a weak version of Sirakov's conjecture in  \cite{CZ1}, and  proved that for all $\beta$ satisfying  \eqref{eq1.06}, $(w_1(x),w_2(x))$ is the unique positive {\em ground state} of \eqref{sys1} up to a translation.


For the singularly perturbed  problem   \eqref{eq1.1}, the existing results  are mainly on the existence and concentration of solutions as $\eps\to0^+$, i.e. searching for the so called peak solutions  in the natural energy space $\mathcal{X}_\eps$. Here, $\vec u=(u_1,u_2)\in \mathcal{X}_\eps $ if
\begin{equation*}
 \Big \{u_i\in  H^1(\R^N):\ \int _{\R^N}  V_i( x)|u_i(x)|^2\dx<\infty \Big\}, \ i=1,2.
\end{equation*}
The norm  of $\vec u=(u_1,u_2)\in \mathcal{X}_\eps $ is defined by
\begin{equation*}
\|(u_1,u_2)\|_{\mathcal{X}_\eps}^2=\sum_{i=1}^2\inte\eps^2|\nabla u_i|^2+V_i( x)|u_i|^2dx.
\end{equation*}
 $\vec u(x)\in \mathcal{X}_\eps$  is called  a solution of \eqref{eq1.1}, provided it is a critical point of the following energy functional in $\mathcal{X}_\eps$
\begin{equation}\label{eq1.7}
J_\eps(\vec u):=\frac{1}{2}\|(u_1,u_2)\|_{\mathcal{X}_\eps}^2-\frac{1}{4}\inte \big(a_1|u_1|^4+a_2|u_2|^4 +2\beta |u_1|^2|u_2|^2\big) \dx.
\end{equation}

Under suitable assumptions  on the limit at infinity of  $V_i(x)$, Lin and Wei \cite{LW2}   proved the existence of positive ground state  for the case of $\beta\in(-\infty, \beta_0)$, where $\beta_0$ is a small unknown positive  constant. Pomponio \cite{Pom} generalized the above results to the case of $a_i=a_i(x)(i=1,2,)$ and $\beta<0$.
In \cite{MPS}, Montefusco, Pellacci and  Squassina considered the case of $a_1=a_2=1$ and $\beta>0$ under the hypothesis that
\begin{equation}\label{cond1}
 0<\min_{B(z,r)} V_i(x)< \min_{\partial B(z,r)}V_i(x),
\text{ for some } z\in \mathbb{R}^3 \text{ and } r>0, \ i=1,2.
\end{equation}
By applying the  mountain pass lemma combining with the technique of truncating the nonlinear term, they   proved that  equation \eqref{eq1.1} possesses a nonnegative solution $(u_{1\eps}, u_{2\eps})$ and the sum $u_{1\eps}+ u_{2\eps}$ has a unique maximum provided $\eps>0$ is small.  Moreover, they showed that if $\beta>0$ is suitably small, then one of $u_{i\eps}$ vanishes as $\eps\to0^+$.  Assuming $V_i(x)$ satisfies an abstract assumption which is a weaker version of \eqref{cond1}, Ikoma and Tanaka \cite{IT} obtained the existence of positive solutions of \eqref{eq1.1}  if $\beta>0$ is suitably small.  By employing the Hardy inequality in $\R^3$ and some truncating technique, Chen and Zou \cite{CZ3}  extended the  work of \cite{IT} to  the case of decaying potentials and $\beta>0$ is relatively large.  Wang and Shi \cite{WS} proposed  the following universe assumption on $V_i(x)$,
\begin{equation}
0<\inf V_i(x)<\liminf_{|x|\to\infty} V_i(x), \ \{y\in\R^N: V_i(y)=\inf_{x\in\R^N} V_i(x),\text{for }i=1,2\}\not=\emptyset,
\end{equation}
and obtained a positive ground state of \eqref{eq1.1} for the case of $\beta>0$ large via the Nehari method. They also proved the existence of multiple nontrivial solutions via the Ljusternik-Schinirelmann category theory. For more results on the existence, multiplicity and concentration of solutions for \eqref{eq1.1} in  bounded domains or whole $\R^N$, one can infer \cite{BC,DW,LW,MMP,Royo,GZZ2,GZZ3,ZZZ} and the references therein.

We note that the solutions obtained by \cite{LW2,MPS,IT,WS} are generally single peak solutions, whose precise definition can be stated  as following.
  \begin{defn}
We call $\big(u_{1\eps}(x),u_{2\eps}(x)\big)\in \mathcal{X}_\eps$ a positive single peak solution of \eqref{eq1.1} concentrating on some point $x_0\in\R^N$, which means that $u_{i\eps}(x)>0(i=1,2,)$ and  $u_\eps(x):=u_{1\eps}(x)+u_{2\eps}(x)$ has a unique maximal point $x_\eps\to x_0$ as  $\eps\to0^+$, while $u_{\eps}(x)\to 0$ uniformly in $\R^N\setminus B_{\theta}(x_0)$ for any $\theta>0$.
\end{defn}

 One natural question is that  whether these single peak solutions are unique under certain assumptions on potentials. As far as we know, there is few result on this respect. As a consequence, in this manuscript we intend to investigate the uniqueness of peak solutions for some typical potentials provided that $\beta>\max\{a_1,a_2\}$. For this purpose, we assume that

 \begin{itemize}
   \item [$(V)$ ] $V_i(x)\in C^\alpha(\R^N; \R^+)$, $\inf_{x\in\R^N}V_i(x)>0(i=1,2,)$ and there exit $\mu>0$ and $z\in\R^3$, such that $V_1(z)=V_2(z)=\mu$.
 \end{itemize}
 Our first result can be stated as follows.
\begin{thm}\label{Thm0}
 Let $\beta>\max\{a_1,a_2\}$ and $(V)$ be satisfied. Assume that   there exist  $m_i,r_0>0$ and $p_i>1$, such that
\begin{equation}\label{eq1.5}
V_i(x)=\mu+m_i |x-z|^{p_i}+o(|x-z|^{p_i}) \ \text{ for all }  x\in B_{r_0}(z), i=1,2,
\end{equation}
and there exists $q_i> p_i-1$ such that
\begin{equation}\label{eq1.6}
\frac{\partial V_i(x+z)}{\partial x_j}=m_ip_{i}|x|^{p_{i}-2}x_j+R_{ij}(x)\,  \text{ with }\,\  |R_{ij}(x)|\leq C|x|^{q_i} \,  \text{ in }\,    B_{r_0}(0).
\end{equation}
Let $\vec u_{\eps}(x)=(u_{1\eps},u_{2\eps})$ and $\vec v_{\eps}(x)=(v_{1\eps},v_{2\eps})$ be two sequences of  positive single peak solutions of (\ref{eq1.1}) in $\mathcal{X}_\eps$ which both concentrate on $z$ and
\begin{equation}\label{eq1.11}
 \lim\limits_{\eps\to0^+} \vec u_{\eps}(\varepsilon x+x_\eps)= \lim\limits_{\eps\to0^+}\vec v_\eps(\varepsilon x+y_\eps)=\vec w(x)
 \text{ strongly in $H^1(\R^N)\times H^1(\R^N)$},
\end{equation}
where $\vec w(x):=(w_1,w_2)$  with $w_i(x)$  being given by (\ref{sys2}), and $x_\eps$ and $y_\eps$  denote the unique maximum point of $u_{1\eps}+u_{2\eps}$ and $v_{1\eps}+v_{2\eps}$, respectively.
Then, there exists $\eps_0>0$ such that $\vec u_{\eps}(x)\equiv \vec v_{\eps}(x)$ provided that $\eps\in(0, \eps_0)$.
\end{thm}
The proof of Theorem \ref{Thm0} is mainly motivated by the arguments of \cite{Cao,Deng,Grossi,GLW,GLWZ} and the references therein.  We also note that, as mentioned  above, the existence of positive single peak solutions  satisfying \eqref{eq1.11} were derived in different articles. For example, the hypothesis  \eqref{eq1.5} indicates \eqref{cond1}, it then follows from  \cite{MPS} that there exists a positive solution  concentrating on $z$ and satisfying \eqref{eq1.11}  if $a_1=a_2=1$ and $\beta>1$. Moreover, if we assume  additionally that  \begin{equation}\label{eq:V}
\liminf_{|x|\to\infty} V_i(x)> \mu=\inf V_i(x)>0, \ \ i=1,2,
\end{equation}
then the single peak solution was also obtained  by \cite[Theorem 1.2]{WS}, and also by \cite[Theorem 1.3]{IT} provided $\beta>0$ is small. Therefore, Theorem \ref{Thm0} indicates that {\em the single peak solutions obtained in \cite{IT,MPS,WS} are unique} provided that \eqref{eq1.5} and \eqref{eq1.6} are satisfied.

For the reader's convenience,  in what follows  we give  a simple proof for the existence of positive single peak solution (which is indeed ground state) in the first part of Theorem \ref{Thm1}  under the hypothesis \eqref{eq:V}. Moreover, for the so called ring-shaped potentials, i.e., $V_i(x)=V_i(|x|)$ attaints its minimum at some spheres, we  derive that the ground state must concentrate on one flattest common  minimal point of  $V_1(x)$ and $V_2(x)$. In addition, we discuss the uniqueness (up to rotations)  of ground state if more information of  $V_i(x)$ near its minimal points is given.
Potentials of ring-shaped appear in many different areas, for instance,  which  are frequently  used in BEC experiments, see e.g. \cite{Ha,HJ,RA}.  In the sequel, we  assume that  $V_1(x)$ and $V_2(x)$ satisfy \eqref{eq:V} and the common minimal points set  \begin{equation}\label{eqV1}
  \mathcal{Z}:=\{x\in\R^N: V_1(x)=V_2(x)=\mu\}\not=\emptyset
\end{equation}
 consists  of   exactly $l(l\in\mathbb{N}^+)$ common spheres, i.e.,
\begin{equation}\label{eq1.17}
  \mathcal{Z}=\big\{x\in\R^N: |x|=A_j>0, 1\leq j\leq l; A_j\not=A_k,\  \forall \ 1\leq j\not=k\leq l\big\}.
\end{equation}
Moreover, we also assume that
 there exist $r_0>0$, $b_{ij}>0$ and $p_{ij}>1$ such that for $i=1, 2$ and $j=1,2,\cdots,l$,
\begin{equation}\label{eq1.18}
V_i(x)=V_i(|x|)=\mu+b_{ij}\big||x|-A_j\big|^{p_{ij}}+o(\big||x|-A_j\big|^{p_{ij}}) \ \ \text{for any }\big||x|-A_j\big|<r_0.
\end{equation}

We note that any rotation of a single peak solution for \eqref{eq1.1} around the origin is still a  single peak solution  because  $V_i(x)=V_i(|x|)$ is radially symmetric. Therefore, to discuss the uniqueness, we must modulate  the rotations. We note that since $\{x:|x|=A_j>0\}$ is a $N-1$-dimensional manifold in $\R^N$, this causes that  the potential is very degenerate in the tangential space of the sphere $\{x:|x|=A_j>0\}$. As a consequence, one cannot derive the uniqueness via  the arguments of Theorem \ref{Thm0}  for the `step 2' in its proof  is not valid anymore. In this case to obtain  the uniqueness (up to rotations) for  general single  peak solutions becomes very challenge. Indeed, we are even uncertain  whether single peak solution is unique or not. However, for the {\em ground state }  for the  second best, by proceeding some refined energy estimates, we can prove that they preserve  some symmetric properties, which is helpful for proving their uniqueness. For this reason, we define
\begin{equation}\label{eq1.170}
c_\eps:=\inf_{\vec u\in \mathcal{N}_\eps }J_\eps(\vec u),
\end{equation}
where the energy functional  $J_\eps(\vec u)$ is given by \eqref{eq1.7}, and
$\mathcal{N}_\eps$ denotes the Nehari manifold and  is defined as
\begin{equation}\label{eq1:N}
\mathcal{N}_\eps:=\Big\{(u_1,u_2)\in \mathcal{X}_\eps\setminus\{0\}:\|(u_1,u_2)\|_{\mathcal{X}_\eps}^2=\inte \big(a_1|u_1|^4+a_2|u_2|^4 +2\beta |u_1|^2|u_2|^2\big) \dx \Big\}.
\end{equation}
We call $(u_1,u_2)$ a ground state of (\ref{eq1.1}) if
\begin{equation}\label{eq1.13}
(u_1,u_2) \text{ is a solution of \eqref{eq1.1} and }J_\eps(u_1,u_2)=c_\eps.
\end{equation}

Assume that $V_i(x)$ satisfies \eqref{eq1.17} and \eqref{eq1.18}.
Denote $\bar\lambda_j $ $(1\leq j\leq l)$  by
\begin{equation}\label{def:unique.Hy}
\bar\lambda_j =\begin{cases}\mu^{1-\frac{N+p_{1j}}{2}}\gamma_1 b_{1j}\displaystyle
\inte|x_N|^{p_{1j}}w^2(x)dx,\,\  &\text{if }\,\  p_{1j}<p_{2j},\\[3mm]
\mu^{1-\frac{N+p_{1j}}{2}}\big(\gamma_1b_{1j}+\gamma_2b_{2j}\big) \displaystyle\inte |x_N|^{p_{1j}}w^2(x)dx,\,\  &\text{if }\,\  p_{1j}=p_{2j},\\[3mm]
\mu^{1-\frac{N+p_{2j}}{2}}\gamma_2 b_{2j}\displaystyle
\inte|x_N|^{p_{2j}}w^2(x)dx,\,\  &\text{if }\,\  p_{1j}>p_{2j},
\end{cases}
\end{equation}
where $\gamma_i>0 $ is given by \eqref{sys2}.
Let ${p}_j:=\min\big\{{p}_{1j},{p}_{2j}\big\}$ and set
\begin{equation}\label{def:beta.p0}
  {p}_0:=\max_{1\leq j\leq l}{p}_j, \ \Gamma:=\big\{1\leq j\leq l:  p_{j}= p_0\big\}\text{ and }\bar{\lambda}_0:=\min\limits_{ j\in \Gamma }\bar\lambda_j.
\end{equation}
Denote
\begin{equation}\label{def:beta.z0}
\mathcal{Z}_0:=\big\{x\in\R^N: |x|=A_j \text{ with } j\in\Gamma \text{ and } \bar\lambda_j= \bar{\lambda}_0\big\}
\end{equation}
 the set of the flattest common minimum points of $V_1(x)$ and $V_2(x)$.

Under the  above assumptions, we  have the following refined estimates for ground states of (\ref{eq1.1}) as $\eps\to0^+$.

\begin{thm}\label{Thm1}
Suppose that $V_1(x)$ and $V_2(x)$ satisfy \eqref{eq:V}  and \eqref{eqV1}, and $\beta>\max\{a_1,a_2\}$.  Then, for $\eps>0$ is small enough, equation (\ref{eq1.1}) has at least one positive ground state, denoted by $\big(u_{1\eps}(x),u_{2\eps}(x)\big)$, and $u_{1\eps}+u_{2\eps}$  possesses a unique maximum point $x_{\eps}$.
Up to a subsequence, there holds that $\lim\limits_{\eps\to0^+}x_{\eps}=\bar{x}_0\in \mathcal{Z}$ and
\begin{equation}\label{eq:th1}
 \lim\limits_{\eps\to0^+} u_{i{\eps}}(\varepsilon x+x_{\eps})
 =w_i(x) \text{ strongly in $H^1(\R^N)$},\quad i=1,2,
\end{equation}
where $w_i$ is given by (\ref{sys2}).
Moreover, suppose that $V_1(x)$ and $V_2(x)$ satisfy  (\ref{eq1.17}) and (\ref{eq1.18}) additionally,  then
we have $\bar{x}_0\in \mathcal{Z}_0$ and
\begin{equation}\label{lim:beta.V.y0}
  \lim\limits_{\eps\to0^+}\frac{|x_\eps|- |\bar x_0|}{\varepsilon}
  =0,
\end{equation}
and
\begin{equation}
c_\eps=\eps^N\big[\bar c_0+\frac{\bar \lambda_0 }2 \eps^{p_0}+ o(\eps^{p_0})\big].
\end{equation}

\end{thm}

%

Based on Theorem \ref{Thm1}, we finally investigate the uniqueness of ground state  for ring-shaped potentials in the following theorem.
\begin{thm}\label{Thm1.5}
Suppose that $V_1(x)$ and $V_2(x)$ satisfy    (\ref{eq1.17}) and (\ref{eq1.18}).
Let $(u_{1\eps}(x),u_{2\eps}(x))$ be the ground state of (\ref{eq1.1}), and $x_{\eps}$ still denotes the unique maximal point of $u_{1\eps}(x)+u_{2\eps}(x)$.  Then, for $\eps>0$ is small, we have
\begin{itemize}
  \item[\rm(i)]  $u_{1\eps}(x)$ and $u_{2\eps}(x)$ are cylindrically symmetric with respect to the line $\overline{Ox_\eps}$.
  \item[\rm(ii)] suppose $\lim\limits_{\eps\to0^+}x_{\eps}=\bar{x}_0\in \mathcal{Z}_0$ with $|\bar x_0|=A_{j_0}$ for some $1\leq j_0\leq l$, and  there exits $r_0>0$ such that  \begin{equation}\label{eq1.28}
\frac{d V_i(|x|)}{d |x|}=b_{ij_0}p_{ij_0}\Big||x|-A_{j_0}\Big|^{p_{ij_0}-2}(|x|-A_{j_0})+R_i(|x|-A_{j_0}) \text{ for }\big||x|-A_{j_0}\big|<r_0,
\end{equation}
where  $$|R_i(|x|-A_{j_0})|\leq C(\big||x|-A_{j_0}\big|^{\tau_i}) \text{ for some }0<\tau_i<p_{ij_0}-1.$$
Then all ground states of (\ref{eq1.1}) which concentrate on the sphere $\{x\in\R^N: |x|=A_{j_0}\}$  indeed make up the following set
      \begin{equation}
      \big\{(u_{1\eps}(Tx),u_{2\eps}(Tx)): T\in O(N)\big\}
       \end{equation}
\end{itemize}

\end{thm}

\begin{rem}
From above theorem we see that if the flattest  common minimal  set of $V_1(x)$ and $V_2(x)$   contains only one sphere, namely, $\mathcal{Z}_0=\{x\in\R^N: |x|=A_{j_0}\}$, then the ground states of (\ref{eq1.1}) is unique up to rotations. It also deserves to finger  out that in \cite{LPWS}, Luo, etc. investigated the existence  of normalized multiple-peak solutions for the single Schr\"odinger equation, where the potential is also assumed to attain its minimum at a $N-1$-dimensional $C^4$ manifold. Especially, if the potential satisfies some {\em non-degenerate } assumptions near the concentration point and {\em does not } preserve any symmetry in $\R^N$, they verified  the normalized multiple-peak solution is unique.  Unfortunately,  the ring-shaped potential is radially symmetric and does not satisfy the assumptions of \cite{LPWS}. Therefore, the arguments of \cite{LPWS} is not applicable for our case.

\end{rem}


This paper is organized as follows. In Section 2, we shall first establish the local Pohozaev identity in Lemma \ref{lem3.1} for nontrivial solutions of \eqref{eq1.1}. Then, we finish the proof of Theorem \ref{Thm0} in Subsetion 2.1.  In Section 3, we  give  refined concentration of ground states of \eqref{eq1.1} under ring-shaped potentials   and finish the proof of Theorem \ref{Thm1} . In Section 4, we first prove in Lemma \ref{lem4.10} that the ground state of \eqref{eq1.1}  must be cylindrically symmetric provided the potential is  ring-shaped, drive the uniqueness (up to rotations) of ground states, which completes the proof of Theorem \ref{Thm1.5}.

\section{Uniqueness of single peak solutions for polynomial  potentials}\label{Se2}
In this section, we intend to prove Theorem \ref{Thm0}, i.e.,  the uniqueness of single peak solutions of \eqref{eq1.1} under polynomial  potentials.
For this purpose, we first establish the local Pohozaev identity in the following lemma.


\begin{lem}\label{lem3.1}
Let $\vec u_\eps(x)=(u_{1\eps},u_{2\eps})$ be a nontrivial solution of (\ref{eq1.1}), then for any $\Omega\subset\R^N$, we have for $j=1,2,\cdots,N$,
\begin{equation}\label{eq3.100}
\begin{split}
\sum_{i=1}^2\int_{\Omega}\frac{\partial V_i(x)}{\partial x_j}u_{i\eps}^2(x)dx&=\sum_{i=1}^2\int_{\partial\Omega} \Big(V_i(x)u_{i\eps}^2(x)-\frac{a_i}{2}u_{i\eps}^4\Big)\nu_jdS-\beta\int_{\partial\Omega}u_{1\eps}^2u_{2\eps}^2dS\\
&+\eps^2\sum_{i=1}^2\int_{\partial\Omega}\Big(|\nabla u_{i\eps}|^2\nu_j-2\frac{\partial u_{i\eps}}{\partial\nu}\frac{\partial u_{i\eps}}{\partial x_j}\Big)dS,
\end{split}
\end{equation}
where $\nu =(\nu _1,\nu _2,\cdots, \nu_N)$ denotes the outward unit normal of $\partial \Omega$.
\end{lem}

\noindent \textbf{Proof.}
Multiply the first equation of \eqref{eq1.1} by $\frac{\partial  u_{1\eps}}{\partial  x_j}$ and integrate over $\Omega$, we then have
\begin{equation}\arraycolsep=1.5pt\begin{array}{lll}
&&-\eps ^2\displaystyle\int_{\Omega}\frac{\partial  u_{1\eps}}{\partial  x_j}\Delta  u_{1\eps}+\displaystyle\int_{\Omega} V_1(x)\frac{\partial  u_{1\eps}}{\partial  x_j} u_{1\eps}\
=\displaystyle a_1\int_{\Omega} \frac{\partial  u_{1\eps}}{\partial  x_j}   u_{1\eps}^3
+\displaystyle\beta \int_{\Omega} \frac{\partial   u_{1\eps}}{\partial  x_j}   u_{1\eps}  u_{2\eps}^2\\[4mm]
\quad\quad&=&\displaystyle \frac{a_1}{4}\int_{\partial\Omega}   u_{1\eps}^4\nu _jdS+\frac{\beta }{2}\int_{\Omega} \frac{\partial   u_{1\eps}^2}{\partial  x_j}   u_{2\eps}^2.
\end{array}\label{5.2:7}
\end{equation}
Note that
\[\arraycolsep=1.5pt\begin{array}{lll}
 &&-\displaystyle\int_{\Omega}\frac{\partial  u_{1\eps}}{\partial  x_j}\Delta  u_{1\eps}
 =- \displaystyle\int_{\partial\Omega}\frac{\partial  u_{1\eps}}{\partial  x_j}\frac{\partial u_{1\eps}}{\partial  \nu}dS+\displaystyle\frac{1}{2}\int_{\partial\Omega} |\nabla u_{1\eps}|^2\nu _jdS,
\end{array}
\]
and
\[
 \displaystyle\int_{\Omega} V_1(x)\frac{\partial  u_{1\eps}}{\partial  x_j}  u_{1\eps}= \frac{1}{2}\int_{\partial\Omega} V_1(x) u_{1\eps}^2\nu _jdS-\frac{1}{2}\int_{\Omega} \frac{\partial V_1(x)}{\partial  x_j} u_{1\eps}^2.
\]
We then derive from (\ref{5.2:7}) that
\begin{equation}\arraycolsep=1.5pt\begin{array}{lll}
\displaystyle\int_{\Omega} \frac{\partial V_1(x)}{\partial  x_j} u_{1\eps}^2+\displaystyle\beta\int_{\Omega} \frac{\partial  u^2_{1\eps}}{\partial  x_j}  u_{2\eps}^2
&=&-2\eps^2\displaystyle\int_{\partial\Omega}\frac{\partial  u_{1\eps}}{\partial  x_j}\frac{\partial u_{1\eps}}{\partial  \nu}dS+\displaystyle \eps^2\int_{\partial\Omega} |\nabla  u_{1\eps}|^2\nu _jdS \\[4mm]
&& +\displaystyle \int_{\partial\Omega} V_1(x)u_{1\eps}^2\nu _jdS-\displaystyle\frac{ a_1}{2}\int_{\partial\Omega}  u_{1\eps}^4\nu _jdS.
\end{array}\label{5.2:8A}
\end{equation}
Similarly, we derive from the second equation of \eqref{eq1.1} that
\begin{equation}\arraycolsep=1.5pt\begin{array}{lll}
\displaystyle \int_{\Omega} \frac{\partial V_2(x)}{\partial  x_j} u_{2\eps}^2+\displaystyle\beta \int_{\Omega} \frac{\partial  u^2_{2\eps}}{\partial  x_j}  u_{1\eps}^2
&=&-2\eps^2\displaystyle\int_{\partial\Omega}\frac{\partial  u_{2\eps}}{\partial  x_j}\frac{\partial  u_{2\eps}}{\partial  \nu}dS+\displaystyle \eps ^2\int_{\partial\Omega} |\nabla  u_{2\eps}|^2\nu _jdS \\[4mm]
&& +\displaystyle \int_{\partial\Omega} V_2(x) u_{2\eps}^2\nu _jdS
-\displaystyle\frac{ a_2}{2}\int_{\partial\Omega}  u_{2\eps}^4\nu _jdS.
\end{array}\label{5.2:8B}
\end{equation}
Note that
\[\displaystyle\int_{\Omega} \frac{\partial  u^2_{2\eps}}{\partial  x_j} u_{1\eps}^2dx=-\int_{\Omega} \frac{\partial  u^2_{1\eps}}{\partial  x_j}  u_{2\eps}^2dx+\int_{\partial\Omega}    u^2_{2\eps}  u_{1\eps}^2 \nu_jdS.\]
This together with (\ref{5.2:8A}) and (\ref{5.2:8B}) gives (\ref{eq3.100}).\qed

\vskip.1truein

We next derive the following  analytic properties for nonnegative solutions of \eqref{eq1.1}.

\begin{lem}\label{lem4.1}
Suppose that $V_i(x)\in C^1(\R^N)$ satisfies $(V)$ for $i=1,\,2$.
Let $(u_{1\eps},u_{2\eps})$ be a positive solution  of \eqref{eq1.1} concentrating on $z$ and denote the unique maximum point of $u_{1\eps}+u_{2\eps}$ by $x_\eps$. Then, up to subsequence,
\begin{itemize}
  \item [(i)]\begin{equation}\label{lem4.1:1}
 \bar{u}_{i\eps}(x):= u_{i\eps}(\varepsilon x+ x_\eps) \overset{\eps\to0^+}\longrightarrow w_i(x) \text{ uniformly in $\R^N$,}
\end{equation}
where $w_i(x)$ is given by \eqref{sys2}.

  \item [(ii)]There exist $R>0$ and $C>0$  independent of $\eps$ such  that
\begin{equation}\label{2:conexp}
\bar{u}_{i\eps}(x)\le Ce^{-\frac{\mu|x|}{2}} \,\ \ \text{if}\,\ |x|>R,\ i=1,2,
\end{equation}
and
\begin{equation}
|\nabla \bar{u}_{i\eps}(x)|\le Ce^{-\frac{\mu|x|}{2}} \, \text{if}\,\ R<|x|<\frac{R}{\eps}, \ i=1,2.
\label{2:conexp2}
\end{equation}
  \item [(iii)] If $V_i(x) (i=1,2,)$ satisfies (\ref{eq1.6}) additionally, then
  \begin{equation}\label{eq3.9}
  \lim_{\eps\to0^+}{|x_\eps-z|}/{\eps}=0.
  \end{equation}
\end{itemize}

\end{lem}

\noindent\textbf{Proof.}  (\ref{lem4.1:1}) and (\ref{2:conexp}) can be proved similarly as (4.2) and (4.6) in \cite{GLWZ}, so we omit their proof here. In view of $x_\eps\to z$ as $\eps\to0^+$, we deduce from \eqref{2:conexp} that
\begin{equation}\label{eq2.90}\big|V_i\big(\varepsilon x+ x_\eps \big)\bar{u}_{i\eps}(x)\big|\leq C e^{-\frac{\mu}{2}|x|} \,\  \text{for}\,\   R<|x|<\frac{R}{\eps}, \end{equation}
where $C>0$ is independent of $\eps$. Note that $\bar u_{1\eps}$ satisfies
\begin{equation*}
-\Delta \bar{u}_{1\eps} +V_i\big(\varepsilon x+ x_\eps \big)\bar{u}_{1\eps} =a_1\bar{u}_{1\eps}^3 +\beta \bar{u}_{2\eps}^2 \bar{u}_{1\eps}   \,\ \mbox{in}\,\  \R^N.
\end{equation*}
Applying  the  elliptic estimates (see $(3.15)$ in \cite{GT}) to above equation, it then follows from \eqref{2:conexp} and \eqref{eq2.90} that
\begin{equation*}
  |\nabla \bar{u}_{1\eps}(x)|\leq Ce^{-\frac{\mu}{2}|x|}\,\  \text{for}\,\   R<|x|<\frac{R}{\eps}.
\end{equation*}
The estimate for $|\nabla \bar u_{2\eps}|$ can be proved similarly.

We finally prove (\ref{eq3.9}). We only give the proof for the case of $p_1<p_2$, while  the case of $p_1\ge p_2$ can be obtained in a similar way.   Taking $\Omega=B_{r_0}(x_\eps)$, it then follows from Lemma \ref{lem3.1}  that, for any $j=1,2,\cdots,N$,

\begin{equation}\label{eq3.10}
\begin{split}
&\sum_{i=1}^2\int_{B_{r_0}(x_\eps)}\frac{\partial V_i(x)}{\partial x_j}u_{i\eps}^2(x)dx=\sum_{i=1}^2\int_{\partial B_{r_0}(x_\eps)} \Big(V_i(x)u_{i\eps}^2(x)-\frac{a_i}{2}u_{i\eps}^4\Big)\nu_jdS\\
&-\beta\int_{\partial B_{r_0}(x_\eps)}u_{1\eps}^2u_{2\eps}^2dS+\eps^2\sum_{i=1}^2\int_{\partial B_{r_0}(x_\eps)}\Big(|\nabla u_{i\eps}|^2\nu_j-2\frac{\partial u_{i\eps}}{\partial\nu}\frac{\partial u_{i\eps}}{\partial x_j}\Big)dS.
\end{split}
\end{equation}
Note that $x_\eps\overset{\eps\to0^+}\longrightarrow z $, we then derive from  \eqref{2:conexp} and \eqref{2:conexp2} that

\begin{equation}\label{eq3.11}
\begin{split}
&\eps^{1-N}\cdot\text{RHS of (\ref{eq3.10})}\\
&=\sum_{i=1}^2\int_{\partial B_{\frac{r_0}\eps}(0)} \Big(V_i(\eps x+x_\eps)\bar u_{i\eps}^2-\frac{a_i}{2}\bar u_{i\eps}^4\Big)\nu_jdS-\beta\int_{\partial B_{\frac{r_0}\eps}(0)}\bar u_{1\eps}^2\bar u_{2\eps}^2dS\\
&+\sum_{i=1}^2\int_{\partial B_{\frac{r_0}\eps}(0)}\Big(|\nabla \bar u_{i\eps}|^2\nu_j-2\frac{\partial \bar u_{i\eps}}{\partial\nu}\frac{\partial \bar u_{i\eps}}{\partial x_j}\Big)dS\\
&=O(e^{-\frac{\gamma}{2\eps}}).
\end{split}
\end{equation}
On the other hand, it follows from (\ref{eq1.6}) that
\begin{equation}\label{eq3.12}
\begin{split}
&\eps^{-N}\cdot\text{LHS of (\ref{eq3.10})}=\sum_{i=1}^2\int_{B_{\frac{r_0}{\eps}}(0)}\frac{\partial V_i(\eps x+x_\eps)}{\partial x_j}\bar u_{i\eps}^2(x)dx\\
&=\sum_{i=1}^2\int_{B_{\frac{r_0}{\eps}}(0)}\big[m_i p_i|\eps x+x_\eps-z|^{p_i-2}(\eps x+x_\eps-z)_j+R_{ij}(\eps x+x_\eps-z)\big]\bar u_{i\eps}^2(x)dx\\
&=\sum_{i=1}^2\int_{B_{\frac{r_0}{\eps}}(0)}m_i p_i|\eps x+x_\eps-z|^{p_i-2}(\eps x+x_\eps-z)_j\bar u_{i\eps}^2(x)dx+O((|x_\eps-z|+\eps)^{\min\{q_1,q_2\}}).
\end{split}
\end{equation}

We claim that
\begin{equation}\label{eq3.13}
{|x_\eps-z|}/{\eps}\text{ is uniformly bounded as }\eps\to0^+.
\end{equation}
For otherwise, if (\ref{eq3.13}) is false,  we can extract a subsequence such that
\begin{equation}\label{eq3.14}
\lim_{\eps\to0^+}\frac{|x_\eps-z|}{\eps}=+\infty, \text{ and }\lim_{\eps\to0^+}\frac{x_\eps-z}{|x_\eps-z|}= \xi \text{ for some }\xi\in S^{N-1}.
\end{equation}
 Since $p_1<p_2$, it then follows from (\ref{eq3.10})-(\ref{eq3.12}) and (\ref{eq3.14}) that
\begin{equation}\label{eq3.15}
\begin{split}
&\Big|\int_{B_{\frac{r_0}{\eps}}(0)} \big|\frac{\eps}{|x_\eps-z|} x+\frac{x_\eps-z}{|x_\eps-z|}\big|^{p_1-2}\big(\frac{\eps}{|x_\eps-z|} x+\frac{x_\eps-z}{|x_\eps-z|}\big)_j\bar u_{1\eps}^2(x)dx\Big|\\
&=\eps^{-1}|x_\eps-z|^{1-p_1} O(e^{-\frac{\gamma}{2\eps}})+O(|x_\eps-z|^{\min\{q_1,q_2\}-p_1+1})+O(|x_\eps-z|^{p_2-p_1}).
\end{split}
\end{equation}
Letting $\eps\to0^+$, it then follows from (\ref{eq3.14}) and (\ref{eq3.15}) that
\begin{equation*}
\xi_j\inte w_1^2(x)dx=0,\  j=1,2,\cdots,N.
\end{equation*}
This indicate that $\xi=0$, which however contradicts (\ref{eq3.14}). (\ref{eq3.13}) is thus proved.

From  (\ref{eq3.13}), we see that, up to a subsequence, $\lim_{\eps\to0^+}\frac{x_\eps-z}{\eps}=y\in\R^N$. Proceeding  the similar argument as that of (\ref{eq3.15}), we deduce that
\begin{equation*}
\begin{split}
&\Big|\int_{B_{\frac{r_0}{\eps}}(0)}|x+\frac{x_\eps-z}{\eps}|^{p_1-2}( x+\frac{x_\eps-z}{\eps})_j\bar u_{1\eps}^2(x)dx\Big|\\
&=\eps^{-1-p_1} O(e^{-\frac{\gamma}{2\eps}})+O(\eps^{\min\{q_1,q_2\}-p_1+1})+O(\eps^{p_2-p_1}).
\end{split}
\end{equation*}
Letting $\eps\to0^+$, it follows that
$$\inte |x+y|^{p_1-2}(x+y)_jw_1^2(x)dx=0, \  j=1,2,\cdots,N.$$
Noting that $w_1(x)=w_1(|x|)$ is strictly decreasing in $|x|$, one can easily deduce from above that $y=0$, and therefore (\ref{eq3.9}) holds.
\qed

\vskip .1truein


\subsection{Proof of Theorem \ref{Thm0}}

\vskip 0.1truein

This subsection is devoted to the proof of Theorem \ref{Thm0} on the uniqueness of positive solutions for the case of $\beta>\max\{a_1,a_2\}$. We first note from \cite[Lemma 2.2 and Theorem 3.1]{DW} that, the positive solution $(w_1,w_2)$ of \eqref{sys1} is {\em non-degenerate} in the sense that, the solution space of the linearized system for (\ref{sys1}) about $(w_1,w_2)$ satisfying
\begin{equation*} 
\begin{cases}
\mathcal{L}_1(\phi_1,\phi_2):=\Delta \phi_1-  \mu \phi_1 + \displaystyle 3a_1 w_1^2\phi_1 + \displaystyle \beta  w_2^2 \phi_1+\displaystyle 2\beta   w_1 w_2\phi_2&=0   \,\ \mbox{in}\,\  \R^N, \\ 
 \mathcal{L}_2(\phi_2,\phi_1):=\Delta \phi_2-  \mu\phi_2 +\displaystyle 3a_2 w_2^2\phi_2 +  \displaystyle \beta  w_1^2 \phi_2+\displaystyle 2\beta   w_1 w_2\phi_1\, &=0   \,\ \mbox{in}\,\  \R^N,
 \end{cases}
\end{equation*}
is exactly $N$-dimensional, namely,
\begin{equation}\label{uniq:limit-3}
\left(\begin{array}{cc} \phi_1\\[2mm]
 \phi_2\end{array}\right)=\sum _{j=1}^{N}b_j\left(\begin{array}{cc} \frac{\partial w_1}{\partial x_j}\\[2mm] \frac{\partial w_2}{\partial x_j}
\end{array}\right)
\end{equation}
for some constants $b_j$.
Assume that   $(u_{1\eps},u_{2\eps})$ and $(v_{1\eps},v_{2\eps})$ are two different positive solutions of (\ref{eq1.1}) concentrating on $z$. Let $x_{\eps}$ and $y_{\eps}$ be the unique maximum point of $u_{1\eps}+u_{2\eps}$ and $v_{1\eps}+v_{2\eps}$, respectively.
Define
\begin{equation}\label{uniq:a-2}
\bar u_{i\eps}(x):=   u_{i\eps}(\varepsilon_{k} x+x_\eps) \text{ and } \bar v_{i\eps}(x):=   v_{i\eps}(\varepsilon x+x_\eps), \ \  \ i=1,2.
\end{equation}
It then follows from Lemma \ref{lem4.1} that
\begin{equation}\label{uniq:a-3}
\bar u_{i\eps}(x), \bar v_{i\eps}(x)\overset{\eps\to 0^+}\longrightarrow w_i(x):=\sqrt{\gamma_i\mu} w(\sqrt \mu x)   \mbox{ uniformly in $\R^N$,\ }\ \ i=1,2.
\end{equation}
From  (\ref{eq1.1}) and (\ref{uniq:a-2}), we see that $(\bar u_{1\eps},\bar u_{2\eps})$ and $(\bar v_{1\eps},\bar v_{2\eps})$  both satisfy the system
\begin{equation}\label{uniq:a-4}\arraycolsep=1.5pt
\left\{\begin{array}{lll}
- \Delta  u_{1} +V_1(\varepsilon x+x_\eps) u_1 &=\displaystyle a_1 u_1^3 +  \beta u_2^2  u_1   \,\ \mbox{in}\,\  \R^N,\\[2.5mm]
- \Delta  u_2 +V_2(\varepsilon x+x_\eps) u_2  &=\displaystyle a_2 u_2^3 +\displaystyle \beta  u_1^2  u_2 \,  \,\ \mbox{in}\,\  \R^N.\,\
\end{array}\right.
\end{equation}

\begin{lem}\label{lem4.2} There exist $C_1>0$ and $C_2>0$  independent of $\eps$, such that
\begin{equation}\label{step-1:1}
C_1\|\bar v_{2\eps}-\bar u_{2\eps}\| _{L^\infty(\R^N)}\le \| \bar v_{1\eps}-\bar u_{1\eps}\| _{L^\infty(\R^N)}\le C_2 \|\bar v_{2\eps}-\bar u_{2\eps}\| _{L^\infty(\R^N)}  \ \ \mbox{as}\ \ \eps\to0^+.
\end{equation}
\end{lem}

\noindent\textbf{Proof.}
We first prove that the right inequality of (\ref{step-1:1}) holds. On the contrary, suppose that
\begin{equation}\label{step-1:2}
 \liminf _{\eps\to0^+}\frac{ \| \bar v_{1\eps}-\bar u_{1\eps}\| _{L^\infty(\R^N)}}{\|\bar v_{2\eps}-\bar u_{2\eps}\| _{L^\infty(\R^N)}}=+\infty.
\end{equation}
From the first equation of (\ref{uniq:a-4}), we have
\begin{equation}\label{uniq:a-5}\arraycolsep=1.5pt
 \begin{array}{lll}
&- \Delta (\bar v_{1\eps}-\bar u_{1\eps}) +V_1(\varepsilon x+x_\eps)(\bar v_{1\eps}-\bar u_{1\eps})\\[2mm]
 =&\displaystyle a_1(\bar v_{1\eps}^3-\bar u_{1\eps}^3) +\beta \Big[\bar  u_{2\eps}^2( \bar v_{1\eps} -\bar u_{1\eps})+\bar v_{1\eps}( \bar v_{2\eps}^2- \bar u_{2\eps}^2 )\Big]  \,\ \mbox{in}\,\  \R^N.
\end{array}
\end{equation}
Set
\begin{equation}\label{step-1:3}
\bar\zeta_{\eps}(x)=\displaystyle\frac{\bar  v_{1\eps}(x)- \bar  u_{1\eps}(x)}{\|\bar  v_{1\eps}- \bar  u_{1\eps}\|_{L^\infty(\R^N)}}.
\end{equation}
It then yields from (\ref{uniq:a-5}) that
\begin{equation}\label{uniq:a-6c}\arraycolsep=1.5pt
 \begin{array}{lll}
&- \Delta \bar \zeta_{\eps} +V_1(\varepsilon x+x_{\eps})\bar \zeta_{\eps}=\displaystyle a_1\big(\bar v_{1\eps}^2+\bar v_{1\eps}\bar u_{1\eps}+\bar u_{1\eps}^2\big)\bar\zeta_{\eps} \\[2mm]
 &+\displaystyle \beta\Big[\bar  u_{2\eps}^2\bar\zeta_{\eps}+\bar v_{1\eps}\frac{ \bar v_{2\eps}^2- \bar u_{2\eps}^2 }{\|\bar v_{1\eps}-\bar u_{1\eps}\|_{L^\infty(\R^2)}}\Big]  \,\ \mbox{in}\,\  \R^N.
\end{array}
\end{equation}
 Note that $\|\bar \zeta_\eps\|_{L^\infty(\R^N)}\le 1$, one can  apply the standard elliptic regularity theory to derive from (\ref{uniq:a-6c}) and (\ref{step-1:2}) that, there exists  $C>0$  independent of $\eps$, such that  $\|\bar\zeta _\eps\|_{C^{1,\alpha }_{loc}(\R^N)}\le C$ for some $\alp \in (0,1)$. This indicates that, up to subsequence, $\bar\zeta_\eps \overset{\eps\to0^+}\longrightarrow \bar\zeta(x) $ for some $ \bar\zeta (x)$.
From (\ref{uniq:a-6c}) we see that  $\bar \zeta $ satisfies
\begin{equation*}
 \Big[-\Delta  +  \mu-\big(\displaystyle3a_1w_1^2+\beta w_2^2\big)\Big]\bar\zeta(x)=0 \,\ \mbox{in}\,\  \R^N.
\end{equation*}
It then follows  from \eqref{sys2} that
\begin{equation}\label{step-1:7}
 \Big[-\Delta  +1-(3a_1\gamma_1+\beta\gamma_2)w^2(x)\Big]\bar\zeta\big(\frac{x}{\sqrt\mu}\big)=0\mbox{ in }  \R^N.
\end{equation}
Noting that $(3a_1\gamma_1+\beta\gamma_2)\in(1,3)$,
we thus conclude from (\ref{step-1:7}) and \cite[Lemma 2.2]{DW} that
$\bar\zeta \equiv 0$ in $\R^N$.

On the other hand, let $z_\eps$ be a point satisfying $|\bar\zeta_\eps(z_\eps)|=\|\bar\zeta_\eps\|_{L^\infty(\R^N)}=1$. Since both $\bar u_{i\eps}$ and $\bar v_{i\eps}$ decay exponentially as $|x|\to\infty$, it then follows from  (\ref{uniq:a-6c}) and standard elliptic regularity theory  that $|z_\eps|\le C$ uniformly in $\eps$. Consequently, we deduce that  $\bar\zeta_\eps\to \bar\zeta \not\equiv 0$ uniformly on $\R^N$, which  leads to a contradiction. Hence (\ref{step-1:2}) is false and  the inequality  on the right hand side of (\ref{step-1:1}) holds true.

Applying the same argument as above, one can deduce from the second equation of (\ref{uniq:a-4}) that the  inequality on the left hand side of (\ref{step-1:1}) also holds. The proof of Lemma \ref{lem4.2} is completed.
\qed
\vskip .1truein

\begin{lem}\label{lem4.3} Assume that $a_1,a _2\in (0, a^*)$, and let
 \begin{equation}\label{uniq:a-7}
\arraycolsep=1.5pt
 \begin{array}{lll}
  \xi_{i\eps}(x)&=\displaystyle\frac{ v_{i\eps}(x)-  u_{i\eps}(x)}{\| v_{1\eps}-  u_{1\eps}\|^\frac{1}{2}_{L^\infty(\R^2)}\| v_{2\eps}- u_{2\eps}\|^\frac{1}{2}_{L^\infty(\R^2)}}, i=1,2.
\end{array}
\end{equation}
 Then for any fixed $x_0\in\R^N$  and $\bar \delta>0$, there exits $\delta=\delta(\eps)\in(\bar\delta,2\bar\delta)$ such that
\begin{equation}
    \int_{\partial B_\delta (x_0)} \Big( \eps ^2 |\nabla  \xi_{i\eps}|^2+   V_i(x)| \xi_{i\eps}|^2\Big)dS=O( \eps ^N)\,\ \text{as}\,\ \eps\to0^+,\,\   i=1, 2.
\label{5.2:6}
\end{equation}
\end{lem}

\noindent\textbf{Proof.}
We first note from   Lemma \ref{lem4.2} that there exists $C>0$ independent of $\eps$ such that
\begin{equation}\label{eq4.29}
0\le| \xi_{1\eps}(x)|,| \xi_{2\eps}(x)|\leq C<\infty \text{ and }|\xi_{1\eps}(x)  \xi_{2\eps}(x)|\le 1 \text{ in } \R^N.
\end{equation}
Recalling from \eqref{eq1.1}, we see  that $(\xi_{1\eps}, \xi_{2\eps})$ satisfies
\begin{equation*}
 \left\{\begin{array}{lll}
&  -\varepsilon^2\Delta \xi_{1\eps}+V_1(x)\xi_{1\eps}= a_1\big(  v_{1\eps}^2+  v_{1\eps}  u_{1\eps}+ u_{1\eps}^2\big)\xi_{1\eps}
 +\displaystyle\beta\big[   u_{2\eps}^2 \xi_{1\eps}+  v_{1\eps}(   v_{2\eps}+ u_{2\eps}) \xi_{2\eps}\big]\\[4mm]
 & -\varepsilon^2 \Delta  \xi_{2\eps}+V_2(x) \xi_{2\eps}=\displaystyle a_2\big( v_{2\eps}^2+  v_{2\eps}  u_{2\eps}+ u_{2\eps}^2\big) \xi_{2\eps}
 +\beta\big[   u_{1\eps}^2 \xi_{2\eps}+  v_{2\eps}(   v_{1\eps}+  u_{1\eps}) \xi_{1\eps}\big]
\end{array}\right.
\end{equation*}
Multiplying the first equation of  above  by $\xi_{1\eps}$ and integrating over $\R^N$, we  obtain that
\[\arraycolsep=1.5pt\begin{array}{lll}
&&\displaystyle \eps ^2\inte |\nabla  \xi_{1\eps}|^2 +\inte V_1(x)| \xi_{1\eps}|^2 \\[4mm]
&=&\displaystyle a_1\inte \big(  v_{1\eps}^2+  v_{1\eps}  u_{1\eps}+  u_{1\eps}^2\big)| \xi_{1\eps}|^2 + \beta\inte \big[   u_{2\eps}^2| \xi_{1\eps}|^2+  v_{1\eps}(   v_{2\eps}+  u_{2\eps}) \xi_{1\eps}\xi_{2\eps}\big]
\\[4mm]
&\le &\displaystyle  C\inte \big(  v_{1\eps}^2+  v_{1\eps} u_{1\eps}+ u_{1\eps}^2\big)+\displaystyle C \inte \big[  u_{2\eps}^2 +  v_{1\eps}(   v_{2\eps}+  u_{2\eps}) \big]
\\[4mm]
&\leq &\displaystyle C\eps^N\inte \big(\bar  v_{1\eps}^2+\bar v_{1\eps}\bar  u_{1\eps}+\bar  u_{1\eps}^2\big)+\displaystyle C\eps^N\inte \big[\bar  u_{2\eps}^2 +\bar  v_{1\eps}( \bar  v_{2\eps}+\bar  u_{2\eps}) \big]
\\[4mm]
&\le& C\eps ^N\,\ \mbox{as} \,\ \eps\to0^+.
\end{array}\]
In view of \eqref{eq4.29} and the exponentially decay estimate of $\bar u_{i\eps}$ and $\bar v_{i\eps}$ in \eqref{2:conexp},
similar to Lemma A.4 of \cite{LPW}, we can apply the mean value theorem of integrals to derive that, for any $x_0\in\R^N$ and  $\bar \delta>0$, there exits $\delta=\delta(\eps)\in(\bar\delta,2\bar\delta)$ such that
\[
    \int_{\partial B_\delta (x_0)} \Big( \eps ^2 |\nabla  \xi_{1\eps}|^2+  V_1(x)| \xi_{1\eps}|^2\Big)dS\le C\bar\delta^{-1}\eps ^N\,\ \mbox{as} \,\ \eps\to0^+.
\]
This means that (\ref{5.2:6}) holds for $i=1$. One can prove  (\ref{5.2:6}) for $i=2$ in a similar way.
\qed

\vskip.1truein
Based on above lemmas, we are ready to proving Theorem \ref{Thm0}.
 We only give the detailed proof for the case of $N=3$,  since the cases of $N=1,2,$ can be derived similarly.

\noindent \textbf{Proof of Theorem \ref{Thm0}:}
{\em  Step 1.}
Let \begin{equation}\label{eq3.2000}
\bar \xi_{i\eps}(x)=\displaystyle\xi_{i\eps}(\eps x+x_\eps), \ \ i=1,2,
\end{equation}
where $\xi_{i\eps}$ is defined by \eqref{uniq:a-7}.
We claim  that, up to a subsequence,
\[(\bar\xi_{1\eps}, \bar\xi_{2\eps})\to (\bar\xi_{10}, \bar\xi_{20}) \,\ \text{in}\,\  C_{loc}^1(\R^3)  \,\ \text{as}\,\   \eps\to0^+,\]
where $(\bar\xi_{10}, \bar\xi_{20})$ satisfies
\begin{equation}\label{uniq:limit-A3}
\left(\begin{array}{cc}   \bar\xi_{10}\\[2mm]
\bar\xi_{20}\end{array}\right)=\displaystyle \sum _{j=1}^{3}d_j\left(\begin{array}{cc} \frac{\partial w_1}{\partial x_j}\\[2mm] \frac{\partial w_2}{\partial x_j}
\end{array}\right)
\end{equation}
for some constants $d_j$ with $j=1, 2, 3$.

Following \eqref{uniq:a-4}, one can actually check that $(\bar\xi_{1\eps}, \bar\xi_{2\eps})$ satisfies
\begin{equation}\label{step-1:9}
\arraycolsep=1.5pt
 \left\{\begin{array}{lll}
  -\Delta \bar\xi_{1\eps}+V_1(\varepsilon x+x_{\varepsilon})\bar\xi_{1\eps}&=a_1\big(\bar v_{1\eps}^2+\bar v_{1\eps}\bar u_{1\eps}+\bar u_{1\eps}^2\big)\bar\xi_{1\eps} \\[2mm]
 & +\beta\big[\bar  u_{2\eps}^2\bar\xi_{1\eps}+\bar v_{1\eps}( \bar v_{2\eps}+\bar u_{2\eps})\bar\xi_{2\eps}\big]\ \mbox{in}\,\  \R^3,\\[4mm]
  -\Delta \bar\xi_{2\eps}+V_2(\varepsilon x+x_{\eps})\bar\xi_{2\eps}&=a_2\big(\bar v_{2\eps}^2+\bar v_{2\eps}\bar u_{2\eps}+\bar u_{2\eps}^2\big)\bar\xi_{2\eps} \\[2mm]
 &+\beta \big[\bar  u_{1\eps}^2\bar\xi_{2\eps}+\bar v_{2\eps}( \bar v_{1\eps}+\bar u_{1\eps})\bar\xi_{1\eps}\big]\ \mbox{in}\,\  \R^3.
\end{array}\right.
\end{equation}
From \eqref{eq4.29} and \eqref{eq3.2000} we see that
\begin{equation}\label{eq4.35}
\text{$\bar\xi_{1\eps}(x)$ and $\bar\xi_{2\eps}(x)$ are both bounded uniformly in $\R^3$}.
\end{equation}
 The standard elliptic regularity theory then implies  that, there exists $C>0$ independent of $\eps$ such that  $\|\bar\xi _{1\eps}\|_{C^{1,\alpha }_{loc}(\R^3)}\le C$ for some $\alp \in (0,1)$. Therefore, up to a subsequence, it yields from \eqref{step-1:9}  that  $(\bar\xi_{1\eps}, \bar\xi_{2\eps})\overset{\eps\to0^+}\longrightarrow (\bar\xi_{10}, \bar\xi_{20})$ in $C_{loc}^1(\R^3)$, where   $(\bar\xi_{10}, \bar\xi_{20})$ satisfies
\begin{equation*}
\arraycolsep=1.5pt
\left\{\begin{array}{lll}
-\Delta   \bar\xi_{10}+  \mu\bar\xi_{10} - \displaystyle 3a_1w_1^2  \bar\xi_{10} - \beta w_2^2   \bar\xi_{10}- 2\beta  w_1w_2  \bar\xi_{20}&=0   \ \ \mbox{in}\,\  \R^3,\\ [1.8mm]
 -\Delta  \bar\xi_{20}+\mu  \bar\xi_{20} -\displaystyle 3a_2 w_2^2 \bar\xi_{20} -\beta w_1^2   \bar\xi_{20}- 2\beta  w_1 w_2  \bar\xi_{10}&=0   \ \ \mbox{in}\,\  \R^3.
\end{array}\right.
\end{equation*}
This together with (\ref{uniq:limit-3}) indicates that $(\bar\xi_{10}, \bar\xi_{20})$ satisfies (\ref{uniq:limit-A3}) for some constants $d_j$ with $j= 1, 2, 3$.

\vskip 0.05truein

\noindent{\em  Step 2.} The constants $d_1=d_2=d_3=0$ in (\ref{uniq:limit-A3}), i.e., $\bar\xi_{10}=\bar\xi_{20}=0$.

Let $\Omega=B_\delta(x_\eps)$ with $\delta=\delta(\eps)>0$ being given by Lemma \ref{lem4.3}.    Applying Lemma \ref{lem3.1} to $(u_{1\eps},u_{2\eps})$  and $(v_{1\eps},v_{2\eps})$  on $\Omega=B_\delta(x_\eps)$, respectively, one can easily derive  that

\begin{equation}\label{5.2:8}
\displaystyle \intB \Big[\frac{\partial V_1(x)}{\partial  x_j}\big( v_{1\eps}+ u_{1\eps}\big)  \xi _{1\eps}+\frac{\partial V_2(x)}{\partial  x_j}\big(v_{2\eps}+u_{2\eps}\big) \xi _{2\eps}\Big] dx:=\sum_{i=1}^2\mathcal{I}^i_\eps+\mathcal{J}_\eps,
\end{equation}
where we denote
\[\arraycolsep=1.5pt\begin{array}{lll}
\mathcal{I}^i_\eps&=&-2\displaystyle \eps ^2\intPB \Big[\frac{\partial  v_{i\eps}}{\partial  x_j}\frac{\partial  \xi_{i\eps}}{\partial  \nu}+\frac{\partial  \xi_{i\eps}}{\partial  x_j}\frac{\partial  u_{i\eps}}{\partial  \nu}\Big]dS+\eps ^2\displaystyle\intPB\nabla  \xi_{i\eps} \cdot\nabla \big( v_{i\eps}+ u_{i\eps}\big)\nu _jdS \\[4mm]
&&+\displaystyle  \intPB V_i(x)\big( v_{i\eps}+ u_{i\eps}\big)  \xi _{i\eps} \nu _jdS -\displaystyle\frac{ a_i}{2}\intPB \big( v^2_{i\eps}+ u^2_{i\eps}\big) \big( v_{i\eps}+ u_{i\eps}\big) \xi _{i\eps}\nu _jdS,
\end{array}\label{5.2:9A}
\]
and
\[\arraycolsep=1.5pt\begin{array}{lll}
\mathcal{J}_\eps&=&-\displaystyle\beta \intPB   \big[ v^2_{1\eps}  ( v_{2\eps}+ u_{2\eps}) \xi _{2\eps}+ u^2_{2\eps}  ( v_{1\eps}+ u_{1\eps}) \xi _{1\eps}\big]dS.
\end{array}\label{5.2:9B}
\]

In view of the fact that $\nabla \bar v_{i\eps}$ satisfies the exponential decay (\ref{2:conexp2}), we deduce from   Lemma \ref{lem4.3} that
\[\arraycolsep=1.5pt\begin{array}{lll}
&&\displaystyle \eps ^2\intPB \Big|\frac{\partial  v_{i\eps}}{\partial  x_j}\frac{\partial  \xi_{i\eps}}{\partial  \nu}\Big|dS\\[4mm]
&\le &\displaystyle \eps \Big(\intPB \Big|\frac{\partial  v_{1\eps}}{\partial  x_j}\Big|^2dS\Big)^{\frac{1}{2}}\Big(\eps ^2\intPB \Big|\frac{\partial  \xi_{i\eps}}{\partial  \nu}\Big|^2dS\Big)^{\frac{1}{2}}\le C\eps ^{1+\frac{3}{2}}e^{-\frac{C\delta}{\eps}}\,\ \mbox{as} \,\ \eps\to0^+.
\end{array}\label{5.2:9a}
\]
Similarly, employing (\ref{2:conexp}), (\ref{2:conexp2}) and Lemma \ref{lem4.3} again, we can  prove that other terms of $\mathcal{I}^i_\eps$ and $\mathcal{J}_\eps$ can also  be controlled  by the order $o(e^{-\frac{C\delta}{\eps}})$.
Therefore, we conclude from above that
\begin{equation}\label{5.2:9aB}
\sum_{i=1}^2\mathcal{I}^i_\eps+\mathcal{J}_\eps=o(e^{-\frac{C\delta}{\eps}}) \,\ \mbox{as} \,\ \eps\to0^+.
\end{equation}

It then follows from (\ref{eq1.6}),   (\ref{5.2:8}) and (\ref{5.2:9aB}) that,
\begin{align}\label{5.2:10}
 o(e^{-\frac{C\delta}{\eps}})&=\displaystyle \sum_{i=1}^2\intB \frac{\partial V_i(x)}{\partial  x_j}\big( v_{i\eps}+ u_{i\eps}\big)  \xi _{i\eps}\dx\\
 &=\displaystyle\eps^3 \sum_{i=1}^2\int_{B_{\frac{\delta}{\eps}}(0)} \frac{\partial V_i\big(\eps[x+(x_{\eps}-z)/\eps]+z\big)}{\partial  x_j}\big(\bar v_{i\eps}+\bar u_{i\eps}\big) \bar\xi _{i\eps}\dx\nonumber\\
 &=\eps^{3} \sum_{i=1}^2\int_{B_{\frac{\delta}{\eps}}(0)} \bigg\{m_{i}p_{i}\eps^{p_{i}-1}\big|x+\frac{x_{\eps}-z}{\eps}\big|^{p_{i}-2}\big(x+\frac{x_{\eps}-z}{\eps}\big)_j\big(\bar v_{i\eps}+\bar u_{i\eps}\big) \bar\xi _{i\eps}\nonumber\\
 &\quad+R_{ij}\big(\eps x+(x_{\eps}-z)\big)\big(\bar v_{i\eps}+\bar u_{i\eps}\big) \bar\xi _{i\eps}\bigg\} \dx.\label{eq4.56}
  \end{align}
 Moreover, we deduce from  (\ref{eq1.6}),  (\ref{2:conexp}), \eqref{eq3.9} and (\ref{eq4.35})  that
 \begin{equation*}
  \begin{split}
& \Big|\int_{B_{\frac{\delta}{\eps}}(0)}R_{ij}\big(\eps x+(x_{\eps}-z)\big)\big(\bar v_{i\eps}+\bar u_{i\eps}\big) \bar\xi _{i\eps} \dx\Big|\\
&\leq C\eps^{q_i} \int_{B_{\frac{\delta}{\eps}}(0)} \big|x+\frac{x_{\eps}-z}{\eps}\big|^{q_i}\big(\bar v_{i\eps}+\bar u_{i\eps}\big) |\bar\xi _{1\eps}|dx\leq C\eps^{q_i}, \\
 \end{split}
  \end{equation*}
 where $q_i>p_{i}-1$ for $i=1,2$.

When $p_{1}=p_{2}$, it then follows from (\ref{eq4.56})  that
  \begin{equation*}
  \begin{split}
  o(1)= p_1\int_{B_{\frac{\delta}{\eps}}(0)} \big|x+\frac{x_{\eps}-z}{\eps}\big|^{p_{1}-2}&\big(x+\frac{x_{\eps}-z}{\eps}\big)_j\sum_{i=1}^2\Big[m_{i}\big(\bar v_{i\eps}+\bar u_{i\eps}\big) \bar\xi _{i\eps}\Big] \dx.
 \end{split}
  \end{equation*}
On the other hand, we deduce from  \eqref{2:conexp} and (\ref{eq4.35}) that
 \begin{equation*}
  \begin{split}
  o(e^{-\frac{C\delta}{\eps}})= \int_{B^c_{\frac{\delta}{\eps}}(0)} \big|x+\frac{x_{\eps}-z}{\eps}\big|^{p_{1}-2}&\big(x+\frac{x_{\eps}-z}{\eps}\big)_j\sum_{i=1}^2\Big[m_{i}\big(\bar v_{i\eps}+\bar u_{i\eps}\big) \bar\xi _{i\eps}\Big] \dx.
 \end{split}
  \end{equation*}
Therefore,
   \begin{equation}\label{eq4.577}
  \begin{split}
  o(1)= \int_{\R^3} \big|x+\frac{x_{\eps}-z}{\eps}\big|^{p_{1}-2}&\big(x+\frac{x_{\eps}-z}{\eps}\big)_j\sum_{i=1}^2\Big[m_{i}\big(\bar v_{i\eps}+\bar u_{i\eps}\big) \bar\xi _{i\eps}\Big] \dx.
 \end{split}
  \end{equation}
Noting that $w_i(x)=w_i(|x|)$, $\frac{\partial w_i}{\partial x_j}=w^\prime(|x|)\frac{x_j}{|x|}$ and $w'(|x|)<0$.
Setting $\eps\to0^+$, we then derive from  \eqref{uniq:a-3}, \eqref{eq3.9} and (\ref{eq4.577}) that
\[ \arraycolsep=1.5pt\begin{array}{lll}
0&=&\displaystyle\int_{\R^3} \big|x\big|^{p_{1}-2}x_j\Big[m_{1} w_1 \bar\xi _{10}+m_{2}w_2  \bar\xi _{20}\Big] dx\\[3mm]
&=&\displaystyle \int_{\R^3} |x|^{p_1-2}x_j\sum_{i=1}^3d_i(m_{1}w_1\frac{\partial w_1}{\partial x_i}+m_{2}w_2\frac{\partial w_2}{\partial x_i})\dx\\
&=&\displaystyle  \frac{d_j}{3}\int_{\R^3} |x|^{p_1-1}(m_{1}w_1w_1^\prime(|x|)+m_{2}w_2w^\prime_2(|x|)\dx,\quad j=1,\,2,\, 3.
\end{array}\]
This indicates that  $d_1=d_2=d_3=0$. Similarly, if $p_{1}\not =p_{2}$,  we can apply the same arguments as above to derive that $d_1=d_2=d_3=0$. Consequently, we conclude that $\bar \xi_{10}=\bar \xi_{20}=0$.

 \vskip 0.1truein

\noindent{\em  Step 3.} $\bar\xi_{10}=\bar\xi_{20}=0$ cannot occur.
Let $(\bar x_\eps, \bar y_\eps)$   satisfy  $|\bar\xi_{1\eps}(\bar x_\eps)\bar\xi_{2\eps}(\bar y_\eps)|=\|\bar\xi_{1\eps}\bar\xi_{2\eps}\|_{L^\infty(\R^3)}=1$.
 By the exponential decay (\ref{2:conexp}), one can easily deduce from  (\ref{step-1:9})  that $|\bar x_\eps|\le C$ and $|\bar y_\eps|\le C$ uniformly in $\eps$. We thus conclude that $\bar\xi_{i\eps}\to \bar\xi_i\not\equiv 0$ uniformly on $\R^3$ as $\eps\to0^+$.  This however contradicts to the fact that $\bar\xi_{10}=\bar\xi_{20}=0$ on $\R^3$. We thus complete the proof of Theorem \ref{Thm0}.
\qed
\vskip .1truein

\section{Refined Concentration of Ground States for Ring-shaped Potentials}

In this section we intend to study the existence and  limit behavior of ground states of \eqref{eq1.1} under ring-shaped  potentials.
For simplicity of notations, we denote
\begin{equation}
B(u_1,u_2):=\inte \big(a_1|u_1|^4+a_2|u_2|^4 +2\beta |u_1|^2|u_2|^2\big) \dx.
\end{equation}
The Nehari manifold corresponds to equation (\ref{sys1}) is defined as
\begin{equation}
\mathcal{N}^\mu:=\big\{(u_1,u_2)\in H^1(\R^N)\times H^1(\R^N)\backslash(0,0):\sum_{i=1}^2\inte|\nabla u_i|^2+\mu|u_i|^2dx= B(u_1,u_2) \big\}.
\end{equation}
Set
\begin{equation}\label{eq2.4}
 c^\mu:= \frac{1}{4}\inf_{(u_1,u_2)\in\mathcal{N}^\mu}B(u_1,u_2)=\frac{1}{4}\inf_{(u_1,u_2)\in\mathcal{N}^\mu}\sum_{i=1}^2\inte|\nabla u_i|^2+\mu|u_i|^2dx,
\end{equation}
and
\begin{equation}\label{eq2.06}
\bar c_1=\inf_{(u_1,0)\in\mathcal{N}^\mu}\frac{1}{4}B(u_1,0), \ \bar c_2=\inf_{(0,u_2)\in\mathcal{N}^\mu}\frac{1}{4}B(0,u_2).
\end{equation}
Then we have the following lemma, which addresses the relationship between $c^\mu$, $\bar c_1$ and $\bar c_2$.
\begin{lem}\label{lem2.1}
If $\beta>\max\{a_1,a_2\}$, then $c^\mu$ can be attained  by $(w_1, w_2)$, where $(w_1, w_2)$ is the unique positive solution of (\ref{sys1}) and is given by  \eqref{sys2}. Moreover, there holds that
\begin{equation}\label{eq2.08}
c^\mu<\min\{\bar c_1,\bar c_2\},\end{equation}
and
\begin{equation}\label{eq2.5}
c^\mu=\frac{1}{4}B(w_1,w_2)=\frac{\mu^{2-\frac{N}{2}}}{4}(a_1\gamma_1^2+a_2\gamma_2^2+2\beta \gamma_1\gamma_2)\inte|w|^4dx.
\end{equation}

\end{lem}
\noindent{\textbf{Proof.}}
One can easily prove that each $\bar c_i$ can be attained by $\sqrt{\frac{\mu}{a_i}}w(\sqrt{\mu}x) \, (i=1,2,)$ with $w(x)$ being the unique positive solution of (\ref{equ:w}) and
\begin{equation}\label{eq2.07}
\bar c_i=\frac{\mu^{2-\frac{N}{2}}}{4a_i}\inte|w|^4\dx,\quad i=1,2.
\end{equation}
Using  $\beta>\max\{a_1,a_2\}$ and noting that $(w_1, w_2)\in \mathcal{N}^\mu$, one can check that
$$c^\mu\leq \frac{1}{4}B(w_1,w_2)=\frac{\mu^{2-\frac{N}{2}}}{4}(a_1\gamma_1^2+a_2\gamma_2^2+2\beta \gamma_1\gamma_2)\inte|w|^4dx<\min\{\bar c_1,\bar c_2\}.$$
This implies (\ref{eq2.08}). As a consequence, there holds that
$$c^\mu:= \frac{1}{4}\inf_{(u_1,u_2)\in\mathcal{N}^\mu, u_1,u_2\not=0}B(u_1,u_2).$$
It then follows from \cite[Theorem 1]{Si} that (\ref{eq2.5}) holds.\qed

\begin{lem} Let $c_\eps$ be defined by \eqref{eq1.170}. Then,
as $\eps\to 0^+$, there holds that
\begin{equation}\label{eq2.09}
0<c^\mu\leq c_\eps/\eps^N\leq c^\mu+o(1).
\end{equation}
Moreover, if  $V_i(x) \ (i=1,2,)$ satisfies \eqref{eq1.17} and \eqref{eq1.18}, we further have
\begin{equation}\label{eq2.10}
c_\eps/\eps^N\leq c^\mu+\frac{\bar \lambda_0 }2 \eps^{p_0}+ o(\eps^{p_0}),
\end{equation}
where $\bar \lambda_0$ is given by (\ref{def:beta.p0}).
\end{lem}
\noindent \textbf{Proof.}
Choose a cutoff function $0\leq \varphi \in C_0^\infty(\R^2)$ such that $\varphi(x)=1$ for $|x|\leq 1$, and $\varphi(x)=0$ for $|x|\geq 2$. Set
\begin{equation}
\tilde w_{i\eps}(x)=w_i\big(\frac{|x-x_0|}{\eps}\big)\varphi\big(\frac{|x-x_0|}{R}\big),\ \ i=1,2,
\end{equation}
where $x_0\in \mathcal{Z}$ and $R>0$ will be determined later. Recalling that
\begin{equation}\label{ide:w}
  \|w\|_2^2=\|\nabla w\|_2^2=\frac{1}{2} \|w\|_4^4,
\end{equation}
and it follows  from \cite[Proposition 4.1]{GNN} that $w(x)$ decays exponentially
 \begin{equation} \label{decay:w}
w(x) \, , \ |\nabla w(x)| = O(|x|^{-\frac{1}{2}}e^{-|x|}) \,\
\text{as} \,\ |x|\to \infty.
\end{equation}
 Using \eqref{ide:w} and \eqref{decay:w}, direct calculations show that
 \begin{equation}\label{eq2.11}
 \begin{split}
\inte|\nabla \tilde w_{i\eps}|^2dx&=\eps^{N-2}\Big[\inte |\nabla w_i(x)|^2\dx+O(e^{-\frac{\sqrt\lambda R}{\eps}})\Big],\\
\inte| \tilde w_{i\eps}|^4dx&=\eps^{N}\Big[\inte |w_i(x)|^4\dx+O(e^{-\frac{\sqrt\lambda R}{\eps}})\Big],
\end{split}
\end{equation}
and
\begin{equation}\label{eq2.12}
\begin{split}
\inte V_i(x)|\tilde w_{i\eps}|^2dx&=\eps^{N}\inte V_i(\eps x+x_0)| w_i(x)|^2\varphi^2(\eps x/R) \dx\\
&=\eps^N\big(\mu\inte|w_i|^2\dx+o(1)\big).
\end{split}
\end{equation}
Taking $t_\eps>0$ such that $t_\eps(\tilde w_{1\eps}, \tilde w_{2\eps})\in \mathcal{N}_\eps$, where  $\mathcal{N}_\eps$ is given by \eqref{eq1:N}. It then follows from (\ref{eq2.11}) and (\ref{eq2.12}) that  $\lim_{\eps\to0^+}t_\eps=1$. Therefore, we recall from \eqref{sys2} and \eqref{eq2.5} that
\begin{equation}\label{eq2.13}
\begin{split}
c_\eps&\leq J_\eps\big(t_\eps(\tilde w_{1\eps}, \tilde w_{2\eps})\big)=\frac{t_\eps^4}{4}B(\tilde w_{1\eps}, \tilde w_{2\eps})\\
&=\frac{\eps^N}{4}(1+o(1))\Big[B(w_1,w_2)+O(e^{-\frac{\sqrt\lambda R}{\eps}})\Big]\\
&=\big(c^\mu+o(1)\big)\eps^N.
\end{split}
\end{equation}
On the other hand, for any $\delta>0$, taking  $(u_{1\eps},u_{2\eps})\in \mathcal{N}_\eps$ such that
$$\frac{1}{4}\sum_{i=1}^2\inte\eps^2|\nabla u_{i\eps}|^2+V_i(x)|u_{i\eps }|^2\dx\leq  c_\eps+\delta\eps^N.$$
From (\ref{eq:V}), we deduce that there exists $\bar t_\eps\in(0,1]$ such that
$\bar t_\eps(v_{1\eps},v_{2\eps}):=\bar t_\eps\big(u_{1\eps}(\eps x), u_{2\eps}(\eps x)\big)\in \mathcal{N}^\mu$. Thus,
\begin{equation*}
\begin{split}
 c^\mu&\leq \frac{\bar t_\eps^2}{4}\sum_{i=1}^2\inte|\nabla v_{i\eps}|^2+\mu|v_{i\eps}|^2\dx\leq \frac{1}{4}\sum_{i=1}^2\inte|\nabla v_{i\eps}|^2+V_i(\eps x)|v_{i\eps }|^2\dx\\
&=\frac{\eps^{-N}}{4}\sum_{i=1}^2\inte\eps^2|\nabla u_{i\eps}|^2+V_i( x)|u_{i\eps}|^2\dx=\eps^{-N}c_\eps+\delta.
\end{split}
\end{equation*}
Letting $\delta\to0^+$, this together with (\ref{eq2.13}) indicates \eqref{eq2.09}.

If $V_i(x)$ satisfies \eqref{eq1.17} and \eqref{eq1.18}, we choose some $x_0\in \mathcal{Z}_0$, where the set $ \mathcal{Z}_0$ is given by \eqref{def:beta.z0}. Thus, $|x_0|=A_{j_0}$ for some $1\leq j_0\leq l$. Taking $R=\frac{r_0}{2}$ with $r_0>0$ being given by (\ref{eq1.18}), then,
\begin{equation}
\begin{split}
&\inte V_i(x)|\tilde w_{i\eps}|^2dx=\eps^{N}\inte V_i(\eps x+x_0)| w_i(x)|^2\varphi^2(\eps x/R) \dx\\
&=\eps^N\Big[\mu\inte|w_i|^2\dx+b_{ij_0}\eps^{p_{ij_0}}\inte|x_N|^{p_{ij_0}}|w_i|^2\dx+o(\eps^{p_{ij_0}})\Big].
\end{split}
\end{equation}
Furthermore,
\begin{equation}
\begin{split}
t_\eps^2&=\frac{\|(\tilde w_{1\eps}, \tilde w_{2\eps})\|^2_{\mathcal{X}\eps}}{B(\tilde w_{1\eps}, \tilde w_{2\eps})}=\frac{\sum_{i=1}^2\inte(|\nabla w_i|^2+\mu w_i^2)\dx+\bar \lambda_0\eps^{p_0}+o(\eps^{p_0})}{B( w_{1},  w_{2})+O(e^{-\frac{\sqrt\lambda R}{\eps}})}\\
&=1+\frac{\bar \lambda_0\eps^{p_0}}{B( w_{1},  w_{2})}+o(\eps^{p_0}).
\end{split}
\end{equation}
Therefore, we have
\begin{equation}
\begin{split}
c_\eps&\leq J_\eps\big(t_\eps(\tilde w_{1\eps}, \tilde w_{2\eps})\big)=\frac{t_\eps^4}{4}B(\tilde w_{1\eps}, \tilde w_{2\eps})\\
&=\frac{\eps^N}{4}\Big[1+\frac{2\bar \lambda_0\eps^{p_0}}{B( w_{1},  w_{2})}+o(\eps^{p_0})\Big]\Big[B(w_1,w_2)+O(e^{-\frac{\sqrt\lambda R}{\eps}})\Big]\\
&=\big[c^\mu+\frac{\bar \lambda_0}{2} \eps^{p_0}+ o(\eps^{p_0})\big]\eps^N.
\end{split}
\end{equation}
This finishes the proof of \eqref{eq2.10}.
\qed
\vskip.1truein


\begin{lem}\label{lem3.3}
Let $V_1(x)$ and $V_2(x)$ satisfy \eqref{eq:V}  and \eqref{eqV1}, and $\beta>\max\{a_1,a_2\}$. Assume that  (\ref{eq1.1}) has a nonnegative ground state $(u_{1\eps},u_{2\eps})$, and  $u_{i\eps}(x)\not\equiv0$ for both $i=1,2$. Then,
\begin{enumerate}
\item [\rm(i)]
   $u_{1\eps}(x)+u_{2\eps}(x)$ admits at least one  global maximum point, denoted it by $x_\eps$. Furthermore, up to a subsequence, there holds that
    \begin{equation}\label{eq3.19}\text{$x_\eps\to \bar x_0$ for some $\bar x_0\in \mathcal{Z}$ as $\eps\to0^+$.}\end{equation}
\item [\rm(ii)]
Define
\begin{equation}\label{def:beta.w.bar}
  {w}_{i\eps}(x)
 :=u_{i\eps}(\varepsilon x+x_\eps)\geq 0,
\end{equation} then
\begin{equation}\label{eq2.8}
  \lim_{\eps\to0^+}w_{i\eps}(x)= w_i(x) \text{ strongly in } H^1(\R^N), \,\ i=1,2.
\end{equation}
\end{enumerate}
\end{lem}

\noindent {\bf Proof.} \textbf{(i).}
For any fixed $\eps>0$, we deduce from \eqref{eq1.1} that  $u_{i\eps}$ satisfies
\[-\eps^2\Delta u_{i\eps}\leq c_{i\eps}(x) u_{i\eps}\,\ \text{in} \,\ \R^N,\,\   i=1,2,\]
where
$$c_{1\eps}(x)=a_1u_{1\eps}^2+\beta u_{2\eps}^2 \text{ and }c_{2\eps}(x)=a_2u_{2\eps}^2+\beta u_{1\eps}^2.$$
Applying the De Giorgi--Nash--Moser theory (cf. \cite[Theorem 4.1]{HL} or \cite[Theorem 8.15]{GT}), we then see that
$$u_{i\eps}(x)\to 0\,\ \text{as}\,\ |x|\to\infty,  \ i=1, 2.$$
This implies that $u_{1\eps}(x)+u_{2\eps}(x)$ admits at least one  global maximum point, denoted it by $x_\eps$. Let $(w_{1\eps},w_{2\eps})$ be defined by
\eqref{def:beta.w.bar}, then it   satisfies
\begin{equation}\label{eq2.16}
\begin{cases}
 -\Delta w_{1\eps} +V_1(\varepsilon x+x_\varepsilon )w_{1\eps}
 =a_1w_{1\eps} ^3+\beta w_{2\eps} ^2w_{1\eps}    \,\ \mbox{in}\,\  \R^N,\\
 -\Delta w_{2\eps} +V_2(\varepsilon x+x_{\varepsilon})w_{2\eps}
 =a_2w_{2\eps} ^3+\beta w_{1\eps} ^2w_{2\eps}    \,\ \mbox{in}\,\  \R^N.
\end{cases}
\end{equation}
Since $w_{1\eps}(x)+ w_{2\eps}(x)$ attains its maximum at $x=0$, we thus derive from above that
$$\mu \big(w_{1\eps}(0)+ w_{2\eps}(0)\big)\leq a_1w_{1\eps} ^3(0)+a_2w_{2\eps} ^3(0)+\beta w_{1\eps}(0)w_{2\eps} (0)\big(w_{1\eps}(0)+ w_{2\eps}(0)\big).$$
This indicates that there exists $C>0$ independent of $\eps$, such that
\begin{equation*}w_{i\eps}(0)\geq C>0\ \  \text{holds for $i=1$ or $2$.}\end{equation*}
Without loss of generality, we assume that  the above estimate holds as least for $i=1$.  As a consequence, there exists $ 0\not\equiv \bar w_1(x)\geq 0$ and $\bar w_2(x)\geq0$, such that, up to a subsequence,
\begin{equation}\label{eq3.200}
\big(w_{1\eps}(x), w_{2\eps}(x)\big)\overset{\eps}\rightharpoonup \big( \bar w_1(x), \bar w_2(x)\big) \text{ weakly in } H^1(\R^N)\times H^1(\R^N).
\end{equation}
We claim  that $\bar w_2(x)\not\equiv0$. For otherwise if $\bar w_2(x)\equiv0$, it then follows from (\ref{eq:V}) that
\begin{equation}
-\Delta  \bar w_1+\mu \bar w_1\leq a_1\bar  w_1^3 \text{ in  }\R^N.
\end{equation}
This indicates that there exists $t_0\in (0,1]$ such that $t_0( \bar w_1,0)\in \mathcal{N}^\mu$. Therefore,  from (\ref{eq2.06}) and (\ref{eq2.08}) we have
\begin{equation}
\begin{split}
\lim_{\eps\to0^+}c_\eps/\eps^N&=\frac14\lim_{\eps\to0}\sum_{i=1}^2\inte|\nabla w_{i\eps}|^2+V_i(\eps x+x_{\eps})|w_{i\eps}|^2\dx\\
&\geq \frac14\inte|\nabla \bar w_{1}|^2+\mu|\bar w_1|^2\dx\geq \bar c_1>c^\mu.
\end{split}
\end{equation}
This leads to a contradiction in view of (\ref{eq2.09}). Thus $\bar w_2(x)\not\equiv0$, and
\begin{equation}\label{eq3.32}0\leq\bar w_{i}(x)\not\equiv0\ \  \text{ for both $i=1,2$.}\end{equation}

Assume that \eqref{eq3.19} is incorrect, then for any sequence $\{\eps_k\}$,  up to a subsequence, there always holds that
\begin{equation}\label{eq2.26}
\lim_{k\to\infty}|x_{\eps_k}|=\infty \text{ or }\lim_{k\to\infty}x_{\eps_k}=z_0\not\in \mathcal{Z}.
\end{equation}
 But in any case, one can deduce from (\ref{eq:V}) that
 \begin{equation}\label{eq3.28}
\lim_{\eps\to0^+}\sum_{i=1}^2\inte|\nabla w_{i\eps}|^2+V_i(\eps x+x_{\eps})|w_{i\eps}|^2\dx>
\sum_{i=1}^2\inte|\nabla \bar  w_i|^2+\mu|\bar  w_i|^2dx.
\end{equation}
Moreover, from  (\ref{eq2.16}) and (\ref{eq2.26}), we derive that there exists  $t_1\in(0,1]$ such that $t_1\big(\bar w_1, \bar w_2\big)\in \mathcal{N}^\mu$.
It then yields from (\ref{eq2.4}) and \eqref{eq3.28} that
\begin{equation}
\begin{split}
\lim_{\eps\to0}c_\eps/\eps^N&=\frac14\lim_{\eps\to0}\sum_{i=1}^2\inte|\nabla w_{i\eps}|^2+V_i(\eps x+x_{\eps})|w_{i\eps}|^2\dx\\
&> \frac14\sum_{i=1}^2\inte|\bar w_i|^2+\mu| \bar w_i|^2dx\geq c^\mu.
\end{split}
\end{equation}
This contradicts (\ref{eq2.09}). Therefore (i) is proved.

\textbf{(ii).}  From \eqref{eq3.19} and \eqref{eq3.200}, we see that  that $\big(\bar w_1(x),\bar  w_2(x)\big)$ solves (\ref{sys1}). Moreover, it follows from \eqref{eq3.32} and the strong maximum principle that both $ \bar w_i(x)>0$ in $\R^N$. Therefore, by the uniqueness of positive solutions for (\ref{sys1}), we see that $\bar w_i(x)\equiv w_i(x)\ (i=1,2,),$ where $w_i(x)$ is given by \eqref{sys2}.
Hence, it follows from (\ref{eq2.09}) that
\begin{equation}
\begin{split}
c^\mu=\lim_{\eps\to0}c_\eps/\eps^N&=\frac14\lim_{\eps\to0}\sum_{i=1}^2\inte|\nabla w_{i\eps}|^2+V_i(\eps x+\eps y_{\eps})|w_{i\eps}|^2\dx\\
&\geq \frac14\sum_{i=1}^2\inte|\nabla \bar w_i|^2+\mu|\bar w_i|^2dx\geq c^\mu.
\end{split}
\end{equation}
This yields (\ref{eq2.8}), and  the proof of the lemma is complete.\qed
\vskip .1truein

\noindent {\bf Proof of Theorem \ref{Thm1}.}
{\em Step1: Existence of positive ground states.} For any $\eps>0$ fixed,  it yields from Theorem \ref{ThmA} in the appendix that \eqref{eq1.1} has at least one nonnegative ground state $(u_{1\eps},u_{2\eps})\in \mathcal{X}_\eps$.
Set $v_{i\eps}(x):=u_{i\eps} (\eps x)\ (i=1,2,)$,
then  $(v_{1\eps}(x), v_{2\eps}(x))$ satisfies
 \begin{equation}\label{eq1.2}
\begin{cases}
-\Delta v_{1\eps} +V_1(\eps x)v_{1\eps} = a_1v_{1\eps}^3 +\beta v_{2\eps}^2 v_{1\eps}   \,\ \mbox{in}\,\  \R^N,\\
-\Delta v_{2\eps} +V_2(\eps x)v_{2\eps} =a_2v_{2\eps}^3 +\beta v_{1\eps}^2 v_{2\eps}   \,\ \mbox{in}\,\  \R^N.\,\
\end{cases}
\end{equation}
We next prove that $u_{1\eps}(x)\not\equiv0$ and $u_{2\eps}(x)\not\equiv0$ provided that $\eps>0$ is small enough. Indeed, if $u_{2\eps}(x)\equiv0$,   then  $(v_{1\eps},0)$ satisfies (\ref{eq1.2}). From (\ref{eq:V}), we see that there exists $t_\eps\in (0,1]$ such  that $t_\eps(v_{1\eps},0)\in \mathcal{N}^\mu$.  As a consequence of Lemma \ref{lem2.1}, we further have
\begin{equation}
\begin{split}
\eps^{-N}c_\eps&=\frac{\eps^{-N}}{4}\inte|\nabla u_{1\eps}|^2+V(x)|u_{1\eps}|^2dx=\frac{1}{4}\inte|\nabla v_{1\eps}|^2+V(\eps x)|v_{1\eps}|^2dx\\
&\geq \frac{1}{4}\inte|\nabla v_{1\eps}|^2+\mu|v_{1\eps}|^2dx\geq \bar c_1>c^\mu.
\end{split}
\end{equation}
This however contradicts (\ref{eq2.09}). Thus, $u_{2\eps}(x)\not\equiv0$. Similarly, we also deduce that $u_{1\eps}(x)\not\equiv0$. Furthermore, by applying the strong maximum principle, we see that $u_{1\eps}(x)>0$ and $u_{2\eps}(x)>0$ in whole $\R^N$.

\vskip .1truein

{\em Step2: Concentration of ground states.}
Let $(w_{1\eps}(x),w_{2\eps}(x))$ be defined by \eqref{def:beta.w.bar}, then \eqref{eq:th1} follows from Lemma \ref{lem3.3}. It remains to prove the uniqueness of the maximum point of $x_\eps$. Indeed, since $V_i(x)$ is H\"older continuous in $\R^N$, it then follows from \eqref{eq3.19} and the standard elliptic regularity theory that
\begin{equation}\label{lim:beta.wbar.C}
{w}_{i{\eps}}(x) \overset{\eps\to0^+}\longrightarrow {w}_{i}(x)\,\ \text{in}\,\  C^2_{loc}(\R^N)\,\ \ \text{for}\,\  i=1, 2.
\end{equation}
Since the origin $x=0$ is the unique maximum point of ${w}_{1}+{w}_{2}$,
 proceeding similar arguments  of \cite[Theorem 1.3]{GLWZ}, one can prove that $x=0$ is the unique maximum point of
 ${w}_{1\eps}+{w}_{2\eps}$, which indicates that $x_\eps$ is the unique maximum point of $u_{1\eps}(x)+u_{2\eps}(x)$.

\vskip 0.1 truein

{\em Step3: Refined energy estimate for ring-shaped potentials.}  Let $(u_{1\eps},u_{2\eps})$ be the convergent subsequence, and $(w_{1\eps},w_{2\eps})$ be given by \eqref{def:beta.w.bar}. Choose  $t_\eps>0$ such that $t_\eps(w_{1\eps},w_{2\eps})\in \mathcal{N}^\mu$, i.e.,
\begin{equation}\label{eq2.46}
t_\eps^2\sum_{i=1}^2\inte|\nabla  w_{i\eps}|^2+\mu|w_{i\eps}|^2\dx=B\big(t_\eps( w_{1\eps}, w_{2\eps})\big).\end{equation}
It then follows from \eqref{eq:V} and \eqref{eq2.16} that
\begin{equation*}
\begin{split}
t_\eps^2&=\frac{\sum_{i=1}^2\inte|\nabla  w_{i\eps}|^2+\mu| w_{i\eps}|^2\dx}{B( w_{1\eps}, w_{2\eps})}
\in(0,1).
\end{split}
\end{equation*}
Applying  \eqref{eq2.16} again, we have
\begin{align}
  c_{\eps}=&\frac{1}{4}\|(u_{1\eps},u_{2\eps})\|_{\mathcal{X}_\eps}^2=\frac{\eps^N}{4}\sum_{i=1}^2\inte|\nabla w_{i\eps}|^2+V_i(\varepsilon x+x_\eps)|{w}_{i\eps}(x)|^2\dx\nonumber\\
  =&\frac{\eps^N}{4}\Big[\sum_{i=1}^2\inte(|\nabla  w_{i\eps}|^2+\mu | w_{i\eps}|^2)+ (V_i(\varepsilon x+x_\eps)-\mu)|{w}_{i\eps}(x)|^2\dx\Big]\label{eq2.41}\\
  \geq &\frac{\eps^N}{4}\Big[\sum_{i=1}^2\inte t_\eps^2(|\nabla  w_{i\eps}|^2+\mu | w_{i\eps}|^2)+ (V_i(\varepsilon x+x_\eps)-\mu)|{w}_{i\eps}(x)|^2\dx\Big]\nonumber\\
  \geq&  c^\mu\eps^N+\frac{\eps^N}{4}\sum_{i=1}^2\inte (V_i(\varepsilon x+x_\eps)-\mu)|{w}_{i\eps}(x)|^2\dx,\label{sub:beta.e.exp}
\end{align}
where the fact $t_\eps( w_{1\eps},  w_{2\eps})\in \mathcal{N}^\mu$ is used in the last inequality.

Since $\lim\limits_{\eps\to0^+}x_{\eps}=\bar{x}_0\in \mathcal{Z}$, where $\mathcal{Z}$ is defined in \eqref{eqV1}, there exists  $1\leq j_0\leq l$ such that
  $|\bar{x}_0|=A_{j_0}$. Without loss of generality,  suppose that ${p}_{j_0}={p}_{1j_0}\leq {p}_{2j_0}$.
We now claim that
\begin{equation}\label{eq3.35}
\text{$\frac{|x_\eps|- A_{j_0}}{\varepsilon}$ \,\,\,is uniformly bounded as $\eps\to0^+$,}\end{equation}
and \begin{equation} \label{eq3.36}  {p}_{j_0}= {p}_0 \text{ where $p_0$ is defined by \eqref{def:beta.p0}}.
\end{equation}
Indeed, since  $\big||\eps x+x_\eps|-A_{j_0}\big|\leq r_0$ for any $x\in B_{\frac{r_0}{2\eps}}(0)$ provided  $\eps>0$ is small enough, we deduce from (\ref{eq1.18})  that
\begin{equation}\label{eq3.430}\begin{split}
&\int_{\R^N}  (V_1(\varepsilon x+x_\eps)-\mu)|{w}_{1\eps}(x)|^2\dx\geq\frac{b_{1j_0}}2\int_{B_{\frac{r_0}{2\eps}}(0)}\big||\eps x+x_\eps|-A_{j_0}\big|^{p_{j_0}}|{w}_{1\eps}(x)|^2\dx\\
&=\frac{b_{1j_0}\eps^{p_{j_0}}}{2}\int_{B_{\frac{r_0}{2\eps}}(0)}\big|\frac{\eps x^2+2x\cdot x_\eps+(x_\eps^2-A_{j_0}^2)/\eps}{|\eps x+x_\eps|+A_{j_0}}\big|^{p_{j_0}}|{w}_{1\eps}(x)|^2\dx.
\end{split}\end{equation}
On the contrary,  suppose that either  $p_{j_0}<{p}_0$ or \eqref{eq3.35} does not hold.  Then, for any given $M>0$, we obtain  from (\ref{eq3.430}) that
\begin{equation*}
\begin{split}
 & \liminf\limits_{\eps\to0^+}\varepsilon^{-{p}_{0}}\int_{\R^N}  (V_1(\varepsilon x+x_\eps)-\mu)|{w}_{1\eps}(x)|^2\dx\\
\geq\quad &
\liminf_{\eps\to0^+}\frac{b_{1j_0}\eps^{p_{j_0}-p_0}}{2}\int_{B_{\frac{r_0}{2\eps}}(0)}\big|\frac{\eps x^2+2x\cdot x_\eps+(x_\eps^2-A_{j_0}^2)/\eps}{|\eps x+x_\eps|+A_{j_0}}\big|^{p_{j_0}}|{w}_{1\eps}(x)|^2\dx
 \geq  M.
\end{split}
\end{equation*}
This together with  \eqref{sub:beta.e.exp} gives that
\begin{equation*}\label{sub:c.beta.e}
 c_{\eps}
   \geq\varepsilon^{N} \big(c^\mu+ M \varepsilon^{{p}_{0}}\big),
\end{equation*}
which however contradicts \eqref{eq2.10}. Thus, the claims   \eqref{eq3.35} and \eqref{eq3.36} are proved.

Following  \eqref{eq3.35}, up to a subsequence of $\{\eps\}$, there exists  $\kappa\in\R$ such that
\begin{equation}\begin{split}
\label{sub:c.beta.ze}
\lim\limits_{\eps\to0^+}\frac{|x_\eps|-|A_{j_0}|}{\varepsilon}=\kappa.
\end{split}\end{equation}
Without loss of generality, we assume that $\bar x_0=(0,0,\cdots, A_{j_0})$.    Applying the Fatou's Lemma, it then follows  from (\ref{eq1.18}), \eqref{def:unique.Hy}, \eqref{eq3.36} and \eqref{sub:c.beta.ze} that
\begin{equation}\label{sub:beta.V1}
\begin{split}
  &\liminf\limits_{\eps\to0^+}\varepsilon^{-{p}_{0}}
 \sum_{i=1}^2\int_{\R^N}  \big(V_i(\varepsilon x+x_\eps)-\mu\big)|{w}_{i\eps}(x)|^2\dx\\
 &=\liminf\limits_{\eps\to0^+}\sum_{i=1}^2\varepsilon^{p_{ij_0}-{p}_{0}}
\int_{\R^N}  \frac{\big(V_i(\varepsilon x+x_\eps)-\mu\big)}{\big||\eps x+x_\eps|-A_{j_0}\big|^{p_{ij_0}}}\big|\frac{\eps x^2+2x\cdot x_\eps+(x_\eps^2-A_{j_0}^2)/\eps}{|\eps x+x_\eps|+A_{j_0}}\big|^{p_{ij_0}}|{w}_{i\eps}(x)|^2\dx\\
&\geq\liminf\limits_{\eps\to0^+}\sum_{i=1}^2\varepsilon^{p_{ij_0}-{p}_{0}}b_{ij_0}\int_{\R^N}  \big|x_N+\kappa\big|^{p_{ij_0}}|{w}_{i}(x)|^2\dx\\
&\geq\liminf\limits_{\eps\to0^+}\sum_{i=1}^2\varepsilon^{p_{ij_0}-{p}_{0}}b_{ij_0}\int_{\R^N}  \big|x_N\big|^{p_{ij_0}}|{w}_{i}(x)|^2\dx= \bar{\lambda}_{j_0}\geq\bar{\lambda}_{0},
\end{split}
\end{equation}
where  $\bar{\lambda}_{j_0}$ and $ \bar{\lambda}_0$ are defined in (\ref{def:unique.Hy}) and (\ref{def:beta.p0}), respectively.
Here it needs to note that the ``="  hold in last two inequalities of \eqref{sub:beta.V1}   if and only if $\kappa=0$ and $\bar x_{0}\in \mathcal{Z}_0$, accordingly.
Moreover,  from  \eqref{eq:th1}, \eqref{eq2.46} and \eqref{sub:beta.V1}, we  have
\begin{equation*}
\begin{split}
t_\eps^2&=\frac{\sum_{i=1}^2\inte|\nabla  w_{i\eps}|^2+\mu| w_{i\eps}|^2\dx}{B( w_{1\eps}, w_{2\eps})}
=1-\frac{\sum_{i=1}^2\inte \big(V_i(\eps x+x_\eps)-\mu\big)| w_{i\eps}|^2\dx}{B( w_{1}, w_{2})+o(1)}\\
&\leq 1-\frac{\bar \lambda_0\eps^{p_0}+o(\eps^{p_0})}{B( w_{1}, w_{2})+o(1)}.
\end{split}
\end{equation*}
This indicates that
\begin{equation*}
\begin{split}
\frac{1}{t_\eps^2}\geq  1+\frac{\bar \lambda_0\eps^{p_0}+o(\eps^{p_0})}{B( w_{1}, w_{2})+o(1)}\big(1+o(1)\big).
\end{split}
\end{equation*}
Noting  that $t_\eps( w_{1\eps}, w_{2\eps})\in \mathcal{N}^\mu$,  it then follows from (\ref{eq2.4}) that
\begin{equation*}
\sum_{i=1}^2\inte(|\nabla  w_{i\eps}|^2+\mu | w_{i\eps}|^2\dx\geq \frac{4 c^\mu}{t_\eps^2}=4c^\mu+\bar \lambda_0\eps^{p_0}+o(\eps^{p_0}).
\end{equation*}
Together with \eqref{eq2.41} and \eqref{sub:beta.V1}, this yields that
\begin{equation}\label{eq2.50}
c_{\eps}\geq \eps^N\big[c^\mu+\frac{\bar \lambda_0}2 \eps^{p_0}+o(\eps^{p_0})\big].
\end{equation}
Combining  with \eqref{eq2.10}, we deduce that  \eqref{eq2.50} is indeed an identity. Hence, all equalities  in \eqref{sub:beta.V1} hold, which implies that   $\kappa=0$ and $ \bar x_{0}\in \mathcal{Z}_0$.
\eqref{lim:beta.V.y0} then follows from \eqref{sub:c.beta.ze}.\qed

\section{Uniqueness of ground state for ring-shaped potentials}

In this section, we investigate the uniqueness of ground state for ring-shaped potentials. We first introduce the following well known Borsuk-Ulam theorem (see, for example, \cite[Theorem 9]{Spa}).

\begin{lem}[Borsuk-Ulam theorem]\label{lemBU}
For any given continuous function $f:\ S^{n_1}\rightarrow \R^{n_2}$ with $n_1\geq n_2\geq1$, there exists $x\in S^{n_1}$ such that $f(x)=f(-x)$.
\end{lem}

Motivated by the argument of Theorem 1 in \cite{Mar}, we have the following cylindrical symmetry result for ground states, which is the core for the proof of Theorem \ref{Thm1.5}.
\begin{lem}\label{lem4.10}
Suppose that $N=2\text{ or }3$, and $V_i(x)=V_i(|x|)$ is  radially symmetric about the origin.
Let $\vec{u }_\eps(x)=(u_{1\eps}(x),u_{2\eps}(x))$ be one ground state of (\ref{eq1.1}).  Assume that  $u_{1\eps}(x)+u_{2\eps}(x)$ has a unique nonzero maximal point, denoted it by  $x_{\eps}\not=0$.  Then, $u_{1\eps}(x)$ and $u_{2\eps}(x)$ are cylindrically  symmetric with the line $\overline{Ox_\eps}$.

\end{lem}

\noindent \textbf{Proof.} We only give the proof when $N=3$, for the case of $N=2$ can be derived similarly. In view of $V_i(x)=V_i(|x|)$, one see that any rotation of $\vec{u }_\eps(x)$ is still a  ground state of (\ref{eq1.1}). Therefore, without loss of generality, we assume that the point $x_\eps$ lines on the $X_3$-axis.

Let $OX_1X_2$ be the vector space spanned by $X_1$ and $X_2$ axes.  For any  $ e\in S^2\cap OX_1X_2$, taking $ v\in e^\perp \cap S^2$, and setting
$$\Pi_v^:=\{x\in\R^3: x\cdot v=0\},\ \Pi_v^+:=\{x\in\R^3: x\cdot v>0\} \text{ and }\Pi_v^-:=\{x\in\R^3: x\cdot v<0\}.$$
Let $$G(x,\vec u_\eps,\nabla \vec u_\eps):= \sum_{i=1}^2 (\eps^2|\nabla u_{i\eps}|+V_i(x)u_{i\eps}^2)+B(\vec u_\eps),$$
and
$$\varphi(v):=\int_{\Pi_v^+}G(x,\vec u_\eps,\nabla \vec u_\eps)-\int_{\Pi_v^-}G(x,\vec u_\eps,\nabla \vec u_\eps).$$
One can easily see that $$\varphi(v)=-\varphi(-v) \ \forall v\in  e^\perp \cap S^2(\cong S^1).$$
It then follows from Lemma \ref{lemBU} that, there exists $v\in  e^\perp \cap S^2$ such that $\varphi(v)=0$. Thus,
\begin{equation}\label{eq4.1}\int_{\Pi_v^+}G(x,\vec u_\eps,\nabla \vec u_\eps)=\int_{\Pi_v^-}G(x,\vec u_\eps,\nabla \vec u_\eps)=0.\end{equation}
Let
$$\vec u^+_\eps=\begin{cases}\vec u_\eps(x)&\ x\in \Pi_v^+\cup\Pi_v,\\
\vec u(P x)&\ x\in\Pi_v^-,
\end{cases}, \ \  \ \vec u^-_\eps=\begin{cases}\vec u_\eps(x)&\ x\in \Pi_v^-\cup\Pi_v\\
\vec u(P x)&\ x\in\Pi_v^+
\end{cases},
$$
where $P$ denotes the orthogonal projection with respect to the hyperplane $\Pi_v$.
It then follows form (\ref{eq4.1}) that
\begin{equation*}\int_{\R^3}G(x,\vec u^+_\eps,\nabla \vec u^+_\eps)=2\int_{\Pi_v^+}G(x,\vec u_\eps,\nabla \vec u_\eps)=0,\end{equation*}
and \begin{equation*}\int_{\R^3}G(x,\vec u^-_\eps,\nabla \vec u^-_\eps)=2\int_{\Pi_v^-}G(x,\vec u_\eps,\nabla \vec u_\eps)=0.\end{equation*}
This indicates that $\vec u^+_\eps, \vec u^-_\eps\in \mathcal{N}_\eps$, and thus

$$c_\eps=J_\eps(\vec u_\eps)=\frac12\big(J_\eps(\vec u^+_\eps)+ J_\eps(\vec u^-_\eps)\big)\geq c_\eps$$
Therefore, we deduce that $\vec u^+_\eps$ and $\vec u^-_\eps$ are both ground states of (\ref{eq1.1}). Set
$$\vec w_{\eps}(x):=\vec u_\eps-\vec u^+_\eps(x)\text{ and }F(\vec u)=\frac{a_1u_1^4+a_2u_2^4}{4}+\frac{\beta u_1^2u_2^2}{2}.$$
We then derive from above that
$$-\Delta \vec w_{\eps}(x)=\mathbb{V}(x)\vec w_{\eps}(x)+\mathbb{A}(x)\vec w_{\eps}(x),$$
where the matrices $\mathbb{V}(x)=diag\{V_1(x), V_2(x)\}$ and $\mathbb{A}(x)$ is given by

$$\mathbb{A}(x)=\left(
                                           \begin{array}{cc}
                                              \int_0^1\frac{\partial^2 F}{\partial u_1^2}(t\vec u_\eps+(1-t)\vec u^+_\eps)dt & \int_0^1\frac{\partial^2 F}{\partial u_1\partial u_2}(t\vec u_\eps+(1-t)\vec u^+_\eps)dt \\
                                             \int_0^1\frac{\partial^2 F}{\partial u_1\partial u_2}(t\vec u_\eps+(1-t)\vec u^+_\eps)dt& \int_0^1\frac{\partial^2 F}{\partial u_2^2}(t\vec u_\eps+(1-t)\vec u^+_\eps) dt\\
                                           \end{array}
                                         \right).$$
Since $\mathbb{V}(x)\in L_{\rm loc}^\infty(\R^3)$ and $\mathbb{A}(x)\in L^\infty (\R^3)$, and note that $\vec u_\eps(x)\equiv\vec u^+_\eps(x)\ \forall\ x\in \Pi_v^+$, we thus can drive  from the unique continuation principle (see for instance, the appendix in \cite{Lop})  that
$$\vec u_\eps(x)=\vec u^+_\eps(x)\text{ for all}\ x\in \R^3.$$
This indicates that $\vec u_\eps(x)$ is symmetric with respect to the hyperplane $\Pi_v$. Moreover, recalling  from Theorem \ref{Thm1} that $u_{1\eps}(x)+u_{2\eps}(x)$ has a unique maximal point $x_{\eps}\not=0$, we then  see that $x_\eps\in \Pi_v$ and thus $\Pi_v=\overline{Oex_\eps}$, where $\overline{Oex_\eps}$ denotes the  hyperplane consists of the line $\overline{Oe}$ and $\overline{Ox_\eps}$. Therefore, we have obtained that
$$\vec u_\eps(x) \text{ is symmetric with respect to } \overline {Oex_\eps} \text{ for all } e\in S^2\cap OX_1X_2.$$
As a consequence, we know that $\vec u_\eps(x)$  is cylindrically symmetric with respect to the line $\overline{Ox_\eps}$.\qed


\vskip.1truein

Based on Lemma \ref{lem4.10} as well as some techniques carry out in the proof of Theorem \ref{Thm0},  we finally finish the proof of Theorem \ref{Thm1.5}.

 \noindent{\bf Proof of Theorem \ref{Thm1.5}.} We only give the proof for the case of $N=3$. Similar to the proof of Theorem \ref{Thm0}, we assume that   $\vec u(x)=(u_{1\eps},u_{2\eps})$ and $\vec v(x)=(v_{1\eps},v_{2\eps})$ are two different nonnegative ground states   of (\ref{eq1.1}), and let $x_{\eps}$ and $y_{\eps}$ be the unique maximum point of $u_{1\eps}+u_{2\eps}$ and $v_{1\eps}+v_{2\eps}$, respectively. We first assume that $x_\eps$ and $y_\eps$ are both lie on the $X_3$ axis, namely,
 \begin{equation}\label{eq4.2}
 x_{\eps}=(0,0, x_{3\eps}) \text{ and }
y_{\eps}=(0,0, y_{3\eps}).\end{equation}
 From Theorem \ref{Thm1}, we have
\begin{equation}
\lim_{\eps\to0^+}\frac{|x_\eps|-A_{j_0}}{\eps}=\lim_{\eps\to0^+}\frac{|y_\eps|-A_{j_0}}{\eps}=0.
\end{equation}
Let $\bar u_{i\eps}(x)$ and $ \bar v_{i\eps}(x)$ be given by (\ref{uniq:a-2}).
$ \xi_{i\eps}(x)$ and $\bar \xi_{i\eps}(x)$  are still defined  by \eqref{uniq:a-7} and (\ref{eq3.2000}), respectively.  Similar to the  the Step 1 in the proof of Theorem \ref{Thm0}, we know that (\ref{uniq:limit-A3}) still holds. We next prove that $d_3=0$. Repeating the arguments of (\ref{5.2:8}) to (\ref{5.2:10}), we see that
\begin{align}\label{eq4.3}
\displaystyle \sum_{i=1}^2\intB \frac{\partial V_i(x)}{\partial  x_j}\big( v_{i\eps}+ u_{i\eps}\big)  \xi _{i\eps}\dx=o(e^{-\frac{C\delta}{\eps}}).
  \end{align}
From (\ref{eq1.28}), we see that
\begin{align}
& \eps^{-3}\cdot L.H.S. of (\ref{eq4.3})\\
 &=\displaystyle \sum_{i=1}^2\int_{B_{\frac{\delta}{\eps}}(0)} \frac{\partial V_i(\eps x+x_{\eps})}{\partial  x_3}\big(\bar v_{i\eps}+\bar u_{i\eps}\big) \bar\xi _{i\eps}\dx\nonumber\\
 &=\sum_{i=1}^2\int_{B_{\frac{\delta}{\eps}}(0)} \bigg\{p_{ij_0}\eps^{p_{ij_0}-1}\big|\frac{\eps x^2+2x_\eps x+(x_\eps^2-A_{j_0}^2)/\eps}{|\eps x+x_\eps|+A_{j_0}}\big|^{p_{ij_0}-2}\frac{\eps x^2+2x_\eps x+(x_\eps^2-A_{j_0}^2)/\eps}{|\eps x+x_\eps|+A_{j_0}}\nonumber\\
 &\quad\quad\cdot\frac{\eps x_3+x_{3\eps}}{|\eps x+x_\eps|}+R_{i}\big(|\eps x+x_\eps|-A_{j_0}\big)\bigg\}\big(\bar v_{i\eps}+\bar u_{i\eps}\big) \bar\xi _{i\eps} \dx.\label{eq4.5}
  \end{align}
Noting that
\begin{equation}\label{eq4.69}\lim_{\eps\to0^+}\frac{\eps x^2+2x_\eps x+(x_\eps^2-A_{j_0}^2)/\eps}{|\eps x+x_\eps|+A_{j_0}}=x_3 \text{ for a. e. }x\in\R^3,\end{equation}
it then follows from (\ref{eq1.28}) that
\begin{equation}\label{eq4.6}
  \begin{split}
& \Big|\int_{B_{\frac{\delta}{\eps}}(0)}R_{i}\big(\eps x+x_{\eps}-A_{j_0}\big)\big(\bar v_{i\eps}+\bar u_{i\eps}\big) \bar\xi _{i\eps} \dx\Big|\cdot \eps^{-\tau_i}\\
&\leq C \int_{B_{\frac{\delta}{\eps}}(0)} \big|\frac{\eps x^2+2x_\eps x+(x_\eps^2-A_{j_0}^2)/\eps}{|\eps x+x_\eps|+A_{j_0}}\big|^{\tau_i}\big(\bar v_{i\eps}+\bar u_{i\eps}\big) |\bar\xi _{1\eps}|dx\leq C<\infty.
 \end{split}
  \end{equation}
Assume that $p_{1j_0}=p_{2j_0}$, similar to the proof of (\ref{eq4.577}), we further deduce from (\ref{eq4.5})  and (\ref{eq4.6}) that   \begin{equation}\label{eq4.7}
  \begin{split}
  o(1)&= \int_{\R^3} \big|\frac{\eps x^2+2x_\eps x+(x_\eps^2-A_{j_0}^2)/\eps}{|\eps x+x_\eps|+A_{j_0}}\big|^{p_{j_0}-2}\big(\frac{\eps x^2+2x_\eps x+(x_\eps^2-A_{j_0}^2)/\eps}{|\eps x+x_\eps|+A_{j_0}}\big)\\
  &\cdot\frac{\eps x_3+x_{3\eps}}{|\eps x+x_\eps|}\sum_{i=1}^2\big(\bar v_{i\eps}+\bar u_{i\eps}\big) \bar\xi _{i\eps} \dx.
 \end{split}
  \end{equation}
Letting $\eps\to0^+$, we then deduce from (\ref{eq4.69}) and (\ref{eq4.7}) that
\[ \arraycolsep=1.5pt\begin{array}{lll}
0&=&\displaystyle\int_{\R^3} \big|x_3\big|^{p_{j_0}-2}x_3\Big[ w_1 \bar\xi _{10}+w_2  \bar\xi _{20}\Big] dx\\[3mm]
&=&\displaystyle d_3\int_{\R^3} |x|^{p_{j_0}-2}x_3(w_1\frac{\partial w_1}{\partial x_3}+m_{2}w_2\frac{\partial w_2}{\partial x_3})\dx\\
&=&\displaystyle  d_3\int_{\R^3} \frac{|x_3|^{p_{j_0}+1}}{|x|}(w_1w_1^\prime(|x|)+w_2w^\prime_2(|x|)\dx .
\end{array}\]
This indicates that  $d_3=0$. For the case of $p_{1j_0}\not=p_{2j_0}$, proceeding  similar arguments as above, one also have $d_3=0$. Moreover, from (\ref{eq4.2}) and Lemma \ref{lem4.10} we see that  $(\bar\xi_{1\eps}, \bar\xi_{2\eps})$,  and thus $(\bar\xi_{10}, \bar\xi_{20})$ is cylindrically symmetric with the $X_3$ axis. This implies that $d_1=d_2=0$. Therefore, we have $\bar\xi_{10}=\bar\xi_{20}=0$, which  however cannot occur by applying the arguments of step 3 for the proof of Theorem \ref{Thm0}.

From above, we see that the ground state which  is cylindrically symmetric with the $X_3$ axis is unique. We denote it by $\vec u_{\eps}(x)$. Then, for any nonnegative ground state $\vec \omega_{\eps}(x)$, there exists $T_{\omega}\in O(3)$  such that $\vec \omega_{\eps}(T_\omega x)$ is symmetric with the $X_3$ axis.  Note that  $\vec \omega(T_\omega x)$ is still a  ground state of (\ref{eq1.1}), we thus deduce that $\vec u_\eps(x)=\vec \omega(T_\omega x)$ in $\R^3$. The proof of Theorem \ref{Thm1.5} is completed.\qed

\appendix
\section{Appendix}

Motivated by the arguments of \cite{WS,MPS}, we intend to prove the existence of non-negative ground states for \eqref{eq1.1} under general potentials in appendix. The following is our main result.

\begin{thm}\label{ThmA}
Suppose that $V_1(x)$ and $V_2(x)$ satisfy \eqref{eq:V}  and \eqref{eqV1}, and $\beta>\max\{a_1,a_2\}$.  Then, for $\eps>0$ is small enough, equation(\ref{eq1.1}) has at least one nonnegtive ground state.
\end{thm}

Denote $V_i^\infty:=\liminf_{|x|\to\infty} V_i(x)$ ($i=1,2,$) and consider the following problem
\begin{equation}\label{A1.1}
\begin{cases}
-\eps^2\Delta u_{1} +V_1^\infty u_{1} =a_1u_{1}^3 +\beta u_{2}^2 u_{1}   \,\ \mbox{in}\,\  \R^N,\ N\leq 3,\\
-\eps^2\Delta u_{2} +V_2^\infty u_{2} =a_2u_{2}^3 +\beta u_{1}^2 u_{2}   \,\ \mbox{in}\,\  \R^N.
\end{cases}
\end{equation}
The corresponding energy functional is defined by
\begin{equation}
J_\eps^\infty(\vec u):=\frac{1}{2}\sum_{i=1}^2\inte(\eps^2|\nabla u_i|^2+V_i^\infty|u_i|^2)dx-\frac{1}{4}B(u_1,u_2) \dx.
\end{equation}
Set
\begin{equation}
c_\eps^\infty :=\inf_{\vec u\in \mathcal{N}_\eps^\infty}J_\eps^\infty(\vec u), \ \text{ where } \mathcal{N}_\eps^\infty:=\big\{\vec u\not=0: \langle J_\eps^\infty(\vec u),\vec u\rangle=0\big\}.
\end{equation}
Then, we have the following lemma.
\begin{lem}
Let (\ref{eq:V}) be satisfied, we have
\begin{equation}\label{A1.05}
c_\eps^\infty/\eps^N>c^\mu\ \text{for any } \eps>0.
\end{equation}
\end{lem}
\noindent \textbf{Proof.} Consider
\begin{equation}\label{A1.2}
\begin{cases}
-\Delta u_{1} +V_1^\infty u_{1} =a_1u_{1}^3 +\beta u_{2}^2 u_{1}   \,\ \mbox{in}\,\  \R^N,\ N\leq 3,\\
-\Delta u_{2} +V_2^\infty u_{2} =a_2u_{2}^3 +\beta u_{1}^2 u_{2}   \,\ \mbox{in}\,\  \R^N,\,\
\end{cases}
\end{equation}
and define
\begin{equation}
c^\infty :=\inf_{\vec u\in \mathcal{N}^\infty}J^\infty(\vec u), \ \text{ where } \mathcal{N}_\eps^\infty:=\big\{\vec u\not=0: \langle J_\eps^\infty(\vec u),\vec u\rangle=0\big\},
\end{equation}
where $J^\infty(\cdot)$ and $\mathcal{N}^\infty$ denote the energy functional and the Nehari manifold of (\ref{A1.2}), respectively. Via the change of variable $\vec v(x):=\vec u(\eps x)$, one can easily see that
\begin{equation}\label{A1.06}c_\eps^\infty/\eps^N=c^\infty.\end{equation} Therefore, it suffices to prove that
\begin{equation}\label{A1.00}
c^\infty>c^\mu.
\end{equation}
Let $\{\vec u_n\}\subset \mathcal{N}^\infty$ be a minimizing sequence of $c^\infty$, we  then have
\begin{equation}\label{A1.3}
\sum_{i=1}^2\inte(|\nabla u_{in}|^2+V_i^\infty|u_{in}|^2)dx=B(u_{1n},u_{2n}) \overset{n}{\to}4 c^\infty>0.
\end{equation}
From (\ref{A1.3}) one can deduce that there exists $\delta>0$ such that
$$\inte |u_{1n}|^2+|u_{2n}|^2 dx\geq \delta>0.$$
Choosing $t_n>0$ such that $t_n\vec u_n\in \mathcal{N}^\mu$, it then follows from (\ref{A1.3}) that
\begin{equation*}
\begin{split}
t_n^2&=\sum_{i=1}^2\inte(|\nabla u_{in}|^2+V_i^\infty|u_{in}|^2)dx/B(\vec u_n)\\
&=1-\frac{\sum_{i=1}^2(V_i^\infty-\mu)\inte|u_{in}|^2dx}{(4c^\infty+o(1))}\leq 1-\delta_1 \text{ for some }\delta_1>0.
\end{split}
\end{equation*}
As a consequence, we have
\begin{equation*}
\begin{split}
c^\mu\leq J_\mu(t_n\vec u_n)=\frac{t_n^4}{4}B(\vec u_n)\leq (1-\delta_1)^2(c^\infty+0(1))<c^\infty.
\end{split}
\end{equation*}
This implies (\ref{A1.00}) and finishes the proof of this lemma.\qed

\begin{lem}\label{lemA3}
Assume that there exists
\begin{equation}\label{A1.4}
\{\vec u_n\}\subset \mathcal{N}_\eps, \text{ such that }  \ J_\eps(\vec u_n)\overset{n}\to c_\eps, \text{ and } J'_\eps(\vec u_n)\overset{n}\to 0.
\end{equation}
Then, if $\eps>0$ is small enough,
\begin{equation}\label{A1.5}
\text{$\vec u_n\overset{n}\to \vec u_0$ strongly in $\mathcal{X}_\eps$, where  $\vec u_0$ is a ground state of (\ref{eq1.1}).}
\end{equation}
\end{lem}

\noindent \textbf{Proof.} From (\ref{A1.4}) we see that $\{\vec u_n\}$ is bounded in $\mathcal{X}_\eps$. Thus, there exists $\vec u_0\in \mathcal{X}_\eps $ such that  $\vec u_n\overset{n}\rightharpoonup \vec u_0$ weakly in $\mathcal{X}_\eps$.  If $\vec u_0\not=0$, then $\vec u_0\in \mathcal{N}_\eps$, and
$$c_\eps=\lim_{n\to\infty}J_\eps(\vec u_n)=\frac{1}{4}\lim_{n\to\infty}B(\vec u_n)\geq \frac{1}{4}\lim_{n\to\infty}B(\vec u_0)=J_\eps(\vec u_0)\geq c_\eps.$$
This implies (\ref{A1.5}). Therefore, to finish the proof of this lemma it remains to prove that $\vec u_0\not=0$.
On the contrary, if $\vec u_0=0$, then
\begin{equation}\label{A.12}\vec u_n\overset{n}\to 0 \text{ in }L^p_{\rm loc}(\R^N)\times L^p_{\rm loc}(\R^N), \text{ for all } 2\leq p <2^*.
\end{equation}

\textbf{Claim:}
\begin{equation}\label{A.13}
\liminf_{n\to\infty}\sum_{i=1}^2\inte (V_i(x)-V_i^\infty)|u_{in}|^2dx\geq0.
\end{equation}
Actually, for any $\delta>0$, there exists $R(\delta)>0$ such that
$$V_i(x)-V_i^\infty>-\delta \text{ if }|x|\geq R(\delta).$$
This indicates that, there exists $C>0$ such that
$$\int_{|x|\geq R(\delta)} (V_i(x)-V_i^\infty)|u_{in}|^2dx\geq-C\delta.$$
From (\ref{A.12}) we also have
$$ \liminf_{n\to\infty}\sum_{i=1}^2\int_{|x|< R(\delta)} (V_i(x)-V_i^\infty)|u_{in}|^2dx=0.$$
The above two estimates yield the claim (\ref{A.13}) holds.

Noting that $\vec u_n\in \mathcal{N}_\eps$, we thus deduce from (\ref{A.13}) that $\limsup_{n\to\infty}\langle J_\eps^\infty(\vec u_n),\vec u_n\rangle\leq0$.
Hence, there exits $t_n>0$ satisfying  $\limsup_{n\to\infty}t_n\leq1$, such that $t_n\vec u_n\in \mathcal{N}_\eps^\infty$. This implies
$$c_\eps^\infty\leq \lim_{n\to\infty}J_\eps^\infty(t_n\vec u_n)=\frac{1}{4}\lim_{n\to\infty}t_n^4B(\vec u_n)=\frac{1}{4}\lim_{n\to\infty}t_n^4J_\eps(t_n\vec u_n)\leq c_\eps.$$
Together with (\ref{eq2.09}), (\ref{A1.05}) and (\ref{A1.06}), we derive  the following contradiction,
$$c^\mu<c^\infty=\lim_{\eps\to0^+}c^\infty_\eps/\eps^N\leq \lim_{\eps\to0^+}c_\eps/\eps^N=c^\mu.$$
Therefore, there holds  that $\vec u_0\not=0$ and the proof is finished.\qed

\begin{lem}\label{lemA4}
There exists nonnegative $\{\vec u_n\}\subset \mathcal{N}_\eps$ satisfying (\ref{A1.4}).
\end{lem}

\noindent \textbf{Proof.}  Let $$\eta_\eps:=\inf_{\gamma\in\Gamma}\sup_{t\in[0,1]}J_\eps(\gamma(t))$$
be the mountain pass level of (\ref{eq1.1}), where
$$\Gamma:=\big\{\gamma(t)\in C(H^1(\R^N)\times H^1(\R^N)): \gamma(0)=0, J_\eps(\gamma(1))<0\big\}.$$
We claim
\begin{equation}\label{A15}
\eta_\eps=c_\eps.
\end{equation}

On the one hand, for any $\vec u\in \mathcal{N}_\eps$, there exists $t_0>1$ such that $J_\eps(t_0\vec u)<0$. Let $\gamma (t)=t(t_0 \vec u)$, then $\gamma(t)\in \Gamma$ and $J_\eps(\gamma(t))$ attains its maximum at the unique point $t=1/t_0$, i.e., $\max_{t\in[0,1]}  J_\eps(\gamma(t))=J_\eps(\vec u)$. This indicates that
\begin{equation}\label{A16}
c_\eps=\inf_{\vec u\in \mathcal{N}_\eps}J_\eps(\vec u)\geq \eta_\eps.
\end{equation}

On the other hand, it follows from \cite[Theorem B.1]{SZ} that there exists a nonnegative $(PS)_{\eta_\eps}$ sequence $\big\{\vec u_n\geq0\big\}\subset \mathcal{X}_\eps$ such that
$$J'_\eps(\vec u_n)\to 0\ \text{ and } J_\eps(\vec u_n)\to \eta_\eps>0. $$
This indicates that $\big\{\vec u_n\big\}$ is uniformly bounded from below and above in $\mathcal{X}_\eps$, and there exits $t_n\overset{n}\to1$ such that $\big\{t_n\vec u_n\big\}\subset \mathcal{N}_\eps$. Therefore, we have
\begin{equation}\label{A17}
J'_\eps(t_n\vec u_n)=t_nJ'_\eps(\vec u_n)+(t_n-t_n^3)B'(\vec u_n)\overset{n}\to0,\end{equation}
and
\begin{equation}\label{A18}
\eta_\eps=\lim_{n\to\infty}J_\eps(\vec u_n)=\lim_{n\to\infty}J_\eps(t_n\vec u_n)\geq c_\eps.\end{equation}
This together with (\ref{A16}) indicates \eqref{A15}. Moreover, \eqref{A15}, \eqref{A17} and \eqref{A18} complete the proof of this lemma.\qed
\vskip.1truein

\noindent \textbf{Proof of Theorem \ref{ThmA}:} Theorem  \ref{ThmA} follows directly from Lemmas \ref{lemA3} and \ref{lemA4}. \qed

\vskip 0.16truein
\noindent {\bf Acknowledgements:} This research was supported by
    the National Natural Science Foundation of China under grant Nos. 11931012, 11871387 and 12171379.


\begin{thebibliography}{40}
  \bibitem{BC} W. Z. Bao and Y. Y. Cai, {\em Ground states of two-component Bose-Einstein condensates with an internal atomic Josephson junction}, East Asia J. Appl. Math. {\bf 1} (2011), 49--81.


\bibitem{BZZ} T. Bartsch, X. X. Zhong and W. M. Zou, {\em Normalized solutions for a coupled Schr\"odinger system}, Math. Ann., {\bf 380} (2021), 1713--1740.


\bibitem{Cao} D. M. Cao,  S. L. Li and P. Luo, {\em Uniqueness of positive bound states with multi-bump for nonlinear Schr\"odinger equations}, Calc. Var. Partial Differential Equations  {\bf 54} (2015),  no. 4, 4037--4063.

\bibitem{CZ1} Z. J. Chen, W. M. Zou, {\em Positive least energy solutions and phase separation for coupled Schr\"odinger equations
with critical exponent}, Arch. Ration. Mech. Anal.,  {\bf 205} (2012),   515--551.

\bibitem{CZ2} Z. J. Chen, W. M. Zou, {\em An optimal constant for the existence of least energy solutions of a coupled Schr\"odinger
system }, Calc. Var. Partial Differential Equations  {\bf 48} (2013),  no. 3-4, 695--711.

\bibitem{CZ3} Z. J. Chen and W. M. Zou, {\em  Standing waves for coupled  Schr\"{o}dinger equations with decaying potentials}, J. Math. Phys. {\bf 54}, 111505 (2013).



\bibitem{DW} E. N. Dancer and J. C. Wei, {\em Spike solutions in coupled nonlinear Schr\"{o}dinger equations with attractive interaction}, Trans. Amer. Math. Soc.  {\bf 361}  (2009),  no. 3, 1189--1208.

\bibitem{Deng}  Y. B. Deng, C. S. Lin and S. Yan, {\em On the prescribed scalar curvature problem in $\R^N$, local uniqueness and periodicity}, J. Math. Pures Appl.   {\bf 104}  (2015),  no. 6, 1013--1044.



\bibitem{Ha} P. L. Halkyard, {Dynamics in Cold Atomic Gases: Resonant Behaviour of the Quantum Delta-Kicked Accelerator and Bose-Einstein Condensates in Ring Traps}, Ph.D Thesis, Durham University, 2010.
\bibitem{HJ} P. L. Halkyard, M. P. A. Jones, S. A. Gardiner, {Rotational response of two-component Bose-Einstein condensates in ring traps}, Phys. Rev. A {81} (2010) 061602 (R).
\bibitem{GNN}  B. Gidas, W. M. Ni and L. Nirenberg,  {\em Symmetry of positive solutions of nonlinear elliptic equations in $\mathbb{R}^n$}, Mathematical analysis and applications  Part A, Adv. in Math. Suppl. Stud. Vol. {\bf 7}  (1981), 369--402.

\bibitem{GT} D. Gilbarg and  N. S. Trudinger, {\em Elliptic Partial Differential Equations of Second Order}, Springer, (1997).


\bibitem{Grossi} M. Grossi, {\em On the number of single-peak solutions of the nonlinear Schr\"odinger equations}, Ann. Inst H. Poinca\'e Anal. Non Lin\'eaire {\bf 19}  (2002), pp. 261--280.

\bibitem{GLW} Y. J. Guo, C. S. Lin and J. C. Wei, {\em Local uniqueness and refined spike profiles of ground states for two-dimensional attractive Bose-Einstein condensates}, SIAM J. Math. Anal. {\bf 49} (2017), 3671--3715.

\bibitem{GLWZ} Y. J. Guo, S. Li, J. C. Wei and X. Y. Zeng, {\em  Ground states of two-component attractive Bose-Einstein condensates I: Existence and uniqueness }, J. Funct. Anal., {\bf276}(1) (2019), 183--230.



\bibitem{GZZ} Y. J. Guo, X. Y. Zeng and H. S. Zhou, {\em Energy estimates and symmetry breaking  in attractive Bose-Einstein condensates with ring-shaped potentials}, Ann. Inst. H. Poincar\'e Anal. Non Lin\'eaire {\bf 33} (2016), 809--828.

\bibitem{GZZ2} Y. J. Guo, X. Y. Zeng and H. S. Zhou, {\em Blow-up solutions for two coupled Gross-Pitaevskii equations with attractive interactions},  Discrete Contin. Dyn. Syst. A, {\bf 37 (7)} (2017), 3749--3786.

\bibitem{GZZ3} Y. J. Guo, X. Y. Zeng and H. S. Zhou, {\em Blow-up behavior of ground states for a nonlinear Schr\"odinger system with attractive and repulsive interactions}, J. Differential Equations 264 (2018), no. 2, 1411--1441.



\bibitem {HMEWC} D. S. Hall, M. R. Matthews,  J. R. Ensher,  C. E. Wieman and E. A. Cornell, {\em Dynamics of component separation in a binary mixture of Bose-Einstein condensates}, Phys. Rev. Lett. {\bf 81} (1998), 1539--1542.

\bibitem {HL} Q. Han and F. H. Lin, {\em Elliptic Partial Differential Equations, Second Edition}, Courant Lect. Notes Math Vol. 1, Courant Institute of Mathematical Science/AMS, New York, (2011).

\bibitem{IT} N. Ikoma, K. Tanaka, {\em A local mountain pass type result for a system of nonlinear Schr\"odinger equations}, Calc. Var. Partial Differential Equations, {\bf40} (2011) 449--480.

\bibitem {Kwong} M. K. Kwong, {\em Uniqueness of positive solutions of $\Delta u-u+u^p=0$ in $\R^N$}, Arch. Ration. Mech. Anal. {\bf 105} (1989), 243--266.




\bibitem {LWCMP} T. C. Lin and J. C. Wei, {\em  Ground state of N coupled nonlinear Schr\"{o}dinger equations in $\R^N$, $N\le 3$}, Comm. Math. Phys. {\bf 255} (2005), 629--653.  Erratum: Comm. Math. Phys. {\bf 277} (2008), 573--576.


\bibitem {LW} T. C. Lin and J. C. Wei, {\em Spikes in two coupled nonlinear Schr\"{o}dinger equations}, Ann. Inst. H. Poincar\'e Anal. Non Lin\'eaire  {\bf 22} (2005), 403--439.

\bibitem {LW2} T. C. Lin and J. C. Wei, {\em Spikes in two-component systems of nonlinear Schr\"odinger equations with trapping potentials}, J. Diff. Eqns. {\bf  229} (2006), 538--569.

\bibitem{LPW} P. Luo, S. J. Peng and C. H.  Wang, Uniqueness of positive solutions with concentration for the
Schr\"odinger-Newton problem, Calc. Var. Partial Differential Equations, {\bf59} (2020), Paper No. 60, 41 pp.

\bibitem{LPWS}P. Luo, S. J. Peng, J. C. Wei, S. S. Yan, Excited states of Bose-Einstein condensates with degenerate attractive interactions, Calc. Var. Partial Differential Equations, {\bf60} (2021), No. 155.




\bibitem{Lop} O. Lopes,  Radial symmetry of minimizers for some translation and rotation invariant
functionals. J. Differ. Equ. {\bf 124}, (1996), 378--388.



\bibitem{Mar} M. Mari\c{s}, On the symmetry of minimizers, Arch. Rational Mech. Anal., {\bf 192}  (2009), 311--330.






 \bibitem{MMP} L. A. Maia,  E. Montefusco and B. Pellacci, {\em Positive solutions for a weakly coupled nonlinear Schr\"{o}dinger system}, J. Differential Equations, {\bf 229} (2006), 743--767.

  \bibitem{MPS}  E. Montefusco, B. Pellacci and M. Squassina, {\em Semiclassical states for  weakly  coupled nonlinear Schr\"{o}dinger systems}, J. Eur. Math. Soc. {\bf 10} (2008), 47--71.




\bibitem{Pom} A. Pomponio, {\em Coupled nonlinear Schr\"odinger systems with potentials}, J. Differential Equations. {\bf 227} (2006), 258--281.





 \bibitem{Royo} J. Royo-Letelier, {\em Segregation and symmetry breaking of strongly coupled twocomponent Bose-Einstein condensates in a harmonic trap}, Calc. Var. Partial Differential Equations {\bf 49} (2014), 103--124.
\bibitem{RA} C. Ryu, M. F. Andersen, P. Clad\'e, Vasant Natarajan, K. Helmerson, W. D. Phillips, {Observation of persistent flow of a Bose-Einstein condensate in a toroidal trap}, Phys. Rev. Lett. {99} (2007) 260401.
 \bibitem{Si}B. Sirakov, Least energy solitary waves for a system of nonlinear Schr\"odinger equations in $\R^n$, Comm. Math. Phys., {\bf271}(1) (2007), 199--221.

\bibitem{Spa} E. H. Spanier, Algebraic Topology, McGraw-Hill, New York, 1966.

\bibitem{SZ} C. A. Stuart and H. S. Zhou, Existence of guided cylindrical tm-modes in an inhomogeneous
self-focusing dielectric, Math. Models Methods Appl. Sci. {\bf20} (2010), 1681--1719.


 \bibitem{WS} J. Wang, J. P. Shi {\em Standing waves of a weakly coupled Schr\"odinger system with distinct potential functions}, J. Differential Equations {\bf260} (2016) 1830--1864.
\bibitem{WY}J. C. Wei,  W. Yao, Uniqueness of positive solutions to some coupled nonlinear Schr\"odinger equations, Commun. Pure Appl. Anal. {\bf11}(3) (2012), 1003--1011.

 \bibitem{WZZ}J. C. Wei,   X. X. Zhong and W. M. Zou, {\em  On Sirakov's open problem and related topics}  Annali SNS-Pisa,  accepted for publication.


\bibitem{ZZZ} X. Y. Zeng, Y. M. Zhang and H.-S. Zhou, Existence and stability of standing waves for a coupled nonlinear schr\"odinger system, Acta Mathematica Scientia, {\bf35} (2015), 45--70.


 \end{thebibliography}
\end{document}